\documentclass[a4 paper]{article}
\usepackage{fancyhdr}
\usepackage{latexsym}
\usepackage{mathrsfs}
\usepackage{cancel}
\usepackage{color}
\usepackage{multirow}
\usepackage{eso-pic}
\usepackage{dsfont}
\usepackage{amssymb}
\usepackage{graphicx}
\usepackage{setspace}
\usepackage{geometry}
\usepackage{amsthm}
\usepackage{hyperref}
\usepackage[intlimits]{amsmath}
\usepackage{hyperref}
\usepackage[affil-it]{authblk}

\usepackage{amssymb,amsfonts}
\usepackage[all,arc]{xy}
\usepackage{enumerate}
\usepackage{mathrsfs}
\usepackage{tikz-cd}
\usepackage{bbm}

\newcommand{\stackspace}{3}
\newcommand{\rstack}[2][0.5cm]{\;\tikz[baseline, yshift=.65ex]%
    {\foreach \k [evaluate=\k as \r using (.5*#2+.5-\k)*\stackspace] in {1,...,#2}{%
    \ifodd\k{\draw[->](0,\r pt)--(#1,\r pt);}%
    \else{\draw[<-](0,\r pt)--(#1,\r pt);}\fi
    }}\;}
\newcommand{\lstack}[2][0.7cm]{\;\tikz[baseline, yshift=.65ex]%
    {\foreach \k [evaluate=\k as \r using (.5*#2+.5-\k)*\stackspace] in {1,...,#2}{%
    \ifodd\k{\draw[<-](0,\r pt)--(#1,\r pt);}%
    \else{\draw[->](0,\r pt)--(#1,\r pt);}\fi
    }}\;}
    
\makeatletter
\def\cleardoublepage{\clearpage\if@twoside \ifodd\c@page\else
   \hbox{}\thispagestyle{empty}\newpage\addtocounter{page}{-1}
   \if@twocolumn\hbox{}\newpage\fi\fi\fi}
\makeatother

\usepackage{graphicx}
\DeclareSymbolFont{lettersA}{U}{txmia}{m}{it}
\DeclareMathSymbol{\pi}{\mathord}{lettersA}{25}
\DeclareMathSymbol{\muup}{\mathord}{lettersA}{22}

\newcommand{\CO}{{\mathcal O}}
\newcommand{\CK}{{\mathcal K}}
\newcommand{\Gr}{\mathrm {Gr}}
\newcommand{\Grbar}{\overline{\mathrm{Gr}}}
\newcommand{\Coh}{\mathrm{Coh}}

\newcommand{\IndCoh}{\mathrm{IndCoh}}
\newcommand{\CBFN}{\mathcal{R}_{G,V}}

\makeatletter
\newcommand{\colim@}[2]{%
  \vtop{\m@th\ialign{##\cr
    \hfil$#1\operator@font colim$\hfil\cr
    \noalign{\nointerlineskip\kern1.5\ex@}#2\cr
    \noalign{\nointerlineskip\kern-\ex@}\cr}}%
}
\newcommand{\colim}{%
  \mathop{\mathpalette\colim@{\rightarrowfill@\scriptscriptstyle}}\nmlimits@
}
\renewcommand{\varprojlim}{%
  \mathop{\mathpalette\varlim@{\leftarrowfill@\scriptscriptstyle}}\nmlimits@
}
\renewcommand{\varinjlim}{%
  \mathop{\mathpalette\varlim@{\rightarrowfill@\scriptscriptstyle}}\nmlimits@
}
\makeatother

\makeatletter
\newcommand*{\doubleleftarrow}[2]{\mathrel{
  \settowidth{\@tempdima}{$\scriptstyle#1$}
  \settowidth{\@tempdimb}{$\scriptstyle#2$}
  \ifdim\@tempdimb>\@tempdima \@tempdima=\@tempdimb\fi
  \mathop{\vcenter{
    \offinterlineskip\ialign{\hbox to\dimexpr\@tempdima+1em{##}\cr
    \leftarrowfill\cr\noalign{\kern.5ex}
    \leftarrowfill\cr}}}\limits^{\!#1}_{\!#2}}}
\newcommand*{\tripleleftarrow}[1]{\mathrel{
  \settowidth{\@tempdima}{$\scriptstyle#1$}
  \mathop{\vcenter{
    \offinterlineskip\ialign{\hbox to\dimexpr\@tempdima+1em{##}\cr
    \leftarrowfill\cr\noalign{\kern.5ex}
    \leftarrowfill\cr\noalign{\kern.5ex}
    \leftarrowfill\cr}}}\limits^{\!#1}}}
\newcommand*{\triplerightarrow}[1]{\mathrel{
  \settowidth{\@tempdima}{$\scriptstyle#1$}
  \mathop{\vcenter{
    \offinterlineskip\ialign{\hbox to\dimexpr\@tempdima+1em{##}\cr
    \rightarrowfill\cr\noalign{\kern.5ex}
    \rightarrowfill\cr\noalign{\kern.5ex}
    \rightarrowfill\cr}}}\limits^{\!#1}}}
\makeatother

\numberwithin{equation}{section}

\begin{document}
\title{Local Operators of 4d $\mathcal{N}=2$ Gauge Theories From the Affine Grassmannian}
\author{Wenjun Niu}
\affil{Department of Mathematics and Center for Quantum Mathematics and Physics, UC Davis}

\maketitle

\begin{abstract}

We give a new, fully mathematical, construction of the space of local operators in the holomorphic-topological twist of 4d $\mathcal N=2$ gauge theories. It is based on computations of morphism spaces in the DG category of line operators, which (by work of Kapustin-Saulina and Cautis-Williams) may be represented as ind-coherent sheaves on the affine Grassmannian and, more generally, on the $\mathcal R_{G,N}$ spaces of Braverman-Finkelberg-Nakajima. We prove that characters of our morphisms spaces reproduce the Schur indices of 4d $\mathcal N=2$ theories, and that the spaces themselves agree with the 4d $\mathcal N=2$ vertex algebras of Beem-Lemos-Liendo-Peelaers-Rastelli-Van Rees, Oh-Yagi, Butson and Jeong. We also generalize our construction to local operators at junctions of Wilson-'t Hooft lines, and compare the Euler character of the morphism spaces to the Schur indices in the work of Cordova-Gaiotto-Shao.

\end{abstract}

\vspace{2cm}

\tableofcontents

\section{Background}

Given a reductive Lie group $G$ and a representation $V$, one can define a 4d $\mathcal N=2$ supersymmetric gauge theory $T_{G,V}$. Physically, it has gauge group $G_c$ (the compact real form) and hypermultiplet matter in representation $V\oplus V^*$. Such a theory admits a variety of twists, labelled by nilpotent elements in the supersymmetry algebra, which are amenable to mathematical study. The prototypical example is Witten's fully topological twist \cite{witten1988topological}, used in the case $G=SL(2),V=0$ to reproduce Donaldson invariants of 4-manifolds. More generally, Kapustin \cite{kapustin2006holomorphic} introduced a holomorphic-topological (HT) twist, which requires spacetime to locally take the form $\mathbb R^2\times \mathbb C$, and depends topologically on $\mathbb R^2$ and holomorphically on $\mathbb C$. Our goal in this paper is to give a new, fully mathematical definition of the space of \emph{local operators} in the HT twist of 4d $\mathcal N=2$ gauge theories, based on a careful computation in the category of \emph{line operators}.

Let $\text{Ops}_{G,V}$ denote the space of local operators in the HT twist of $T_{G,V}$. It is physically defined as the cohomology of the HT supercharge acting on the full space of local operators in the untwisted theory. It has a $\mathbb Z$-valued cohomological grading `$F$', and an additional non-cohomological $\frac12\mathbb Z$-valued grading `$J$', corresponding to spin in the holomorphic plane (mixed with an $SU(2)$ R-symmetry). This space has been well studied from other perspectives. In particular:
\begin{itemize}
\item Its graded Euler character is the ``Schur index'' of $T_{G,V}$,
\begin{equation} \chi_q\, \text{Ops}_{G,V} := \text{Tr}_{\text{Ops}_{G,V}} (-1)^F q^J \,=\, I_{\text{Schur}}[T_{G,V}]\,. \end{equation}
The Schur index, introduced in \cite{gadde2011four, gadde2013gauge}, is a particular specialization of the 4d $\mathcal N=2$ superconformal index \cite{kinney2007index, romelsberger2006counting}; though it makes sense even when a 4d $\mathcal N=2$ theory is not conformal. In this paper, we will always use $q$ as a formal variable counting the weight of the loop rotation. 

\item The space $\text{Ops}_{G,V}$ itself is the vacuum module of a Poisson vertex algebra $\mathcal V_{G,V}$. This Poisson vertex algebra was constructed for general 4d $\mathcal N=2$ theories from a more physical perspective by Oh and Yagi \cite{oh2020poisson}, and from a mathematical perspective by Dylan Butson \cite{butson2021equivariant}, as a BRST reduction of classical beta-gamma algebras valued in $T^*V$. When $T_{G,V}$ is superconformal --- meaning quadratic indices satisfy $C_2(N)=C_2(\mathfrak g)$ --- the vertex algebra can be further quantized by introducing an Omega background, yielding a VOA $\mathcal V^{\hbar}_{G,V}$. These VOA's were first introduced in superconformal theories by Beem-Lemos-Liendo-Peelaers-Rastelli-Rees \cite{beem2015infinite}. However, the work of \cite{beem2015infinite} did not define this vertex algebra from the point of view of local operators in $T_{G,V}$. The fact that the two pictures coincide is a nontrivial result, and can be explained using a "cigar-like" reduction. This is explained in  \cite{oh2019chiral} and \cite{jeong2019scft} from a physical perspective and  \cite{butson2021equivariant} from a mathematical perspective. This deformation does not alter the underlying vector space of the vacuum module, so
\begin{equation} \text{Ops}_{G,V} \simeq \mathcal V_{G,V} \simeq \mathcal V^\hbar_{G,V}\,. \end{equation}
\end{itemize}

The perspective we will take, similar to  \cite{butson2021equivariant}, uses the category of line operators in the HT twist of a 4d $\mathcal N=2$ theory. 
Physically, the objects of this category are line operators supported on a line in the topological $\mathbb  R^2$ plane and the origin of the holomorphic $\mathbb C$ plane. The category contains half-BPS Wilson-'t Hooft lines, as well as more general quarter-BPS line operators. The category was given a geometric description by Cautis and Williams, in \cite{cautis2019cluster} for pure gauge theory ($V=0$) and \cite{Cautistoappear} for general $V$, following the physical predictions of \cite{kapustin2006holomorphic}.
This category is described as the category of equivariant coherent sheaves on the Braverman-Finkelberg-Nakajima (BFN) space. Let us recall the construction here.

Let $\mathcal{O}$ be the ring of formal power series $\mathbb{C}[\![z]\!]$ and $\mathcal{K}$ the field of formal Laurent series $\mathbb{C}(\!(z)\!)$. Denote by $\mathbb{D}=\mathrm{Spec}(\CO)$ the formal disk and $\mathbb{D}^*=\mathrm{Spec}(\CK)$ the formal punctured disk. Given a Lie group $G$ and representation $V$ of $G$, the BFN space is defined by the following base-change diagram:
\begin{equation}\label{eqBFNdef} 
\begin{tikzcd}
\mathcal{R}_{G,V} \rar \arrow[d] &  V(\mathcal{O})\arrow[d] \\  G(\mathcal{K})\times_{G(\mathcal{O})}V(\mathcal{O}) \rar &  V(\mathcal{K})
\end{tikzcd}
\end{equation} 
This space can be alternatively described as the moduli space of triples $(\mathcal{P},\varphi, s)$ where $\mathcal{P}$ is a principal $G$ torsor over $\mathbb{D}$, $\varphi$ is a trivialization of $\mathcal{P}$ over $\mathbb{D}^*$, and $s$ is a section of the associated $V$ bundle over $\mathbb{D}$ that is sent, under $\varphi$, to a regular section of the trivial $V$ bundle. This remarkable object was first studied in \cite{nakajima2015towards,braverman2016towards}, where the authors used its Borel-Moore homology to construct Coulumb branches of 3d $\mathcal{N}=4$ gauge theories. The K-theoretic version is studied in \cite{braverman2016coulomb}, \cite{cautis2019cluster} as well as \cite{finkelberg2019multiplicative}.  There is an action of $\tilde{G}_\CO$ on the left on $\mathcal{R}_{G,V}$, where $\mathbb{C}^*$ is some cover group of loop rotation. The category of line operators in the HT twist is expected to be the category of equivariant coherent sheaves on $\mathcal{R}_{G,V}$ as a monoidal category. This category is denoted by $\mathrm{Coh}(\tilde{G}_\CO\setminus \CBFN)$, whose monoidal structure is defined carefully in \cite{Cautistoappear}. Under this expectation, the trivial line is the tensor identity $\mathbbm{1}$ of the category and the space of local operators is the endomorphism algebra
\begin{equation} \text{Ops}_{G,V} \simeq \mathrm{End}_{\mathrm{Coh}(\tilde{G}_\CO\setminus \CBFN)}(\mathbbm{1})\,. \end{equation}

This expectation is in fact not easy to work with. Over the last decade, much effort has been made towards using the techniques of derived algebraic geometry to rigorously define aspects of the physical quantum field theories, and this is one example of this. Indeed, when considering sheaves on this space, several issues appear: first, that $\CBFN$ is a non-Noetherian ind-scheme, and so some care needs to be taken when considering coherent sheaves; second, equation \eqref{eqBFNdef} is a base-change diagram, and one should consider this as a diagram of DG ind-schemes; third, that $\CBFN$ is not reduced. If one considers Borel-Moore homology, as in \cite{nakajima2015towards,braverman2016towards}, or K theory, as in \cite{braverman2016coulomb,cautis2019cluster,finkelberg2019multiplicative}, then one can consider the reduced-classical scheme, but this is not enough to properly produce the space of local operators. 

This paper is devoted to such an endeavor, with the help of techniques of derived algebraic geometry. The structure of a DG  indscheme and its category of sheaves are studied in \cite{gaitsgory2014dg, gaitsgory2019study, gaitsgory2020study} for locally almost finite type, and \cite{raskin2020homological} in general. We will use the machinery of \cite{gaitsgory2014dg,gaitsgory2019study, gaitsgory2020study, raskin2020homological} to define the category of ind-coherent sheaves, and compute carefully and explicitly the endomorphism algebra of the identity line up to quasi-isomorphism. We will show (Theorem \ref{Extidentitypure} and Theorem \ref{Extidentitymatter}) that after certain shift of gradings, it coincides with the vacuum module of the Poisson vertex algebra in \cite{oh2020poisson} and \cite{butson2021equivariant}. We will also consider the insertion of other (half-BPS) Wilson-'t Hooft line operators and the space of local operators at their junctions, which will give rise to modules of the algebra $\mathcal{V}_{G,V}$. We will compare the Euler characters of these modules with the index formulae in \cite{cordova2016infrared}. In the category of line operators $\mathrm{Coh}(\tilde{G}_\CO\setminus \CBFN)$, these Wilson-'t Hooft operators are equivariant vector bundles on various $\tilde{G}_\CO$ orbits inside $\CBFN$. 

In the work of \cite{butson2021equivariant}, the author showed that the category of line operators in HT twist has the structure of a factorization $\mathbb{E}_1$-category (\cite[Section 5.11]{butson2021equivariant}). This structure should give rise to the structure of a $\mathbb{E}_2$ factorization algebra to $\text{Ops}_{G,V}$, which is the aforementioned Poisson vertex algebra structure. We hope that this work with the procedure outlined in \cite{butson2021equivariant} can help rigorously produce factorization algebras from the category of line operators. 

The paper is structured in the following way. In Section \ref{2},  we give an introduction to the geometry of the BFN spaces. In Section \ref{3}, we introduce the Poisson vertex algebra $\mathcal{V}_{G,V}$. In Section \ref{4}, we introduce the category of ind-coherent sheaves on the BFN spaces; we compute the space of endomorphisms of the identity line for pure gauge theory (Theorem \ref{Extidentitypure}) and for gauge theory with matter (Theorem \ref{Extidentitymatter}); in both cases the result coincides with the vacuum module of the Poisson vertex algebra~$\mathcal{V}_{G,V}$. In Section \ref{5}, Theorem \ref{Thmfundamental}, we compute the space of local operators at the junction of fundamental t'Hooft lines and dyonic Wilson-'t Hooft lines, in pure $\mathrm{PSL}(2)$ theory, and compare the results with the formulae of Cordova-Gaiotto-Shao~\cite{cordova2016infrared}. 

\vspace{10pt}

\noindent\textbf{Acknowledgements.} I would like to thank Tudor Dimofte for his guidance throughout this research project. I would like to thank Kevin Costello for motivating some of the constructions in this paper (in particular ideas contained in section \ref{AbelianGaugeGroup}), and Justin Hilburn for teaching me about Witt vectors. I would like to thank Harold Williams for very patiently pointing out mistakes in the first draft of the paper, and for teaching me derived algebraic geometry. I would like to thank Dylan Butson for sharing his insights into factorization categories and sharing with me his thesis. I would like to thank Nick Rozenblyum for teaching me about Corollary \ref{CorExtidentity}. I would also like to thank Niklas Garner and Philsang Yoo for many helpful discussions. I would like to thank my friend Don Manuel for many encouragements. 

\section{The Affine Grassmannian and the BFN space}\label{2}

\subsection{The Geometry of the Affine Grassmannian}

We will work with the ground field $\mathbb{C}$. Let $G$ be a reductive Lie group. The affine Grassmannian of $G$ is the quotient:
\begin{equation}
\mathrm{Gr}_G:=\mathcal{R}_{G,0}=G(\mathcal{K})/ G(\mathcal{O}).
\end{equation}
It turns out that $\Gr_G$ is a classical ind-scheme; its geometry is well studied in the literature. We will in this section recall some basic facts about this space. For details, see \cite{zhu2016introduction}. In particular, we note that the study of the geometry of this space has two complications, one is that it is an ind-scheme; the other is that it is not always reduced.

\subsubsection{The Affine Grassmannian of $\mathrm{GL}_n$}

In the case of $\mathrm{GL}_n$, the affine Grassmannian $\Gr_{\mathrm{GL}_n}$ can be defined alternatively as the moduli space of lattices in $\CK^n$. More precisely, if $R$ is an algebra over $\mathbb{C}$,  then an $R$ family of lattices in $\mathcal{K}^n$ is a finitely-generated projective $R[\![z]\!]$-submodule $\Lambda$ of $R(\!(z)\!)^n$ such that $\Lambda\otimes_{R[\![z]\!]}R(\!(z)\!)=R(\!(z)\!)^n$. The affine Grassmannian $\Gr_{\mathrm{GL}_n}$ can be defined as the presheaf assigning to $R$ the set of $R$ families of lattices in $\mathcal{K}^n$. We have:

\newtheorem{IndGrGLn}{Proposition}[section]

\begin{IndGrGLn}\label{IndGrGLn}
$\Gr_{\mathrm{GL}_n}$ is represented by a classical ind-projective ind-scheme. Namely, it can be written as a colimit of classical projective schemes under closed embeddings.
\end{IndGrGLn}

Moreover, $\Gr_{\mathrm{GL}_n}$ is formally smooth in the following sense:

\newtheorem{FormalSmoothDef}[IndGrGLn]{Definition}

\begin{FormalSmoothDef}
An ind-scheme $X=\varinjlim X_n$ is formally smooth if for any algebra $R$ and nilpotent ideal $I\subseteq R$, the map $X(R)\to X(R/I)$ is surjective. 
\end{FormalSmoothDef}

\subsubsection{The Affine Grassmannian of General $G$}

Given a reductive Lie group $G$, choosing a faithful representation $G\to \mathrm{GL}_n$ one obtains a closed embedding $\Gr_G\hookrightarrow \Gr_{\mathrm{GL}_n}$, and from Proposition \ref{IndGrGLn} one concludes that $\Gr_G$ is an ind-projective ind-scheme. There is a canonical isomorphism $\pi_0(\Gr_G)\cong \pi_1(G)$, and the connected components of $\Gr_G$ are labeled by the fundamental group of $G$, which is also the quotient of the co-weight lattice of $G$ by its co-root lattice. All connected components are isomorphic to each other, with an isomorphism given by left multiplication. It turns out, however, this space is not always reduced, as can be seen from the following example:

\vspace{10pt}
\noindent\textbf{Example.}  Consider $T=\mathbb{C}^*$. Then $\Gr_T=\CK^*/\CO^*$. The $\mathbb{C}$ points of this space is a disjoint union of infinitely many copies of $\mathrm{Spec}(\mathbb{C})$. They can be represented by $\{z^n\vert n\in \mathbb{Z}\}$, since any nonzero element in $\CK$ is of the form $z^{n}g[z]$ for some $g[z]\in \CO^*$. However, if we evaluate $\Gr_G$ on the algebra $\mathbb{C}[\epsilon]$ with $\epsilon^2=0$, then an element in $\mathbb{C}[\epsilon]\otimes \CK$ is invertible iff its image under the map $\mathbb{C}[\epsilon]\otimes \CK\to \CK$ setting $\epsilon\mapsto 0$ is invertible. This means that the set of invertible elements in $\mathbb{C}[\epsilon]\otimes \CK$ is $\CK^*\oplus \epsilon \CK$, while the set of invertible elements in $\CO$ is $\CO^*\oplus \epsilon \CO$. The quotient is not a discrete set anymore, but rather an infinite-dimensional vector bundle over $\Gr_T(\mathbb{C})$. The fibre of this bundle at a point $z^n$ is $\CK/z^n\CO$, and should be interpreted as the tangent space of $\Gr_T$ at $z^n$. This presents difficulty in considering the category of coherent sheaves, as in this case, the category of sheaves on $\Gr_T$ is different from the category of sheaves on its $\mathbb{C}$ points, even though the $K_0$ groups of the two categories are isomorphic. 
\vspace{10pt}

In the case when $G$ is semi-simple, $\Gr_G$ is in fact reduced. In general, we denote by $\Gr_{G,red}$ the reduced ind-scheme of $\Gr_G$. In the following, we will present a different stratification of $\Gr_G$ from the one obtained through an embedding $\Gr_G\hookrightarrow \Gr_{\mathrm{GL}_n}$.

\subsubsection{A Stratification for $\Gr_{G,red}$}\label{Stratificationred}

Let $T$ be a maximal torus of $G$ and $B$ a Borel subgroup of $G$ containing $T$. Denote the associated weight lattice by $X^*(T)$ and coweight lattice by $X_*(T)$. The choice of a Borel subgroup determines a set of dominant weights $X^*(T)_+$ and dominant coweights $X_*(T)_+$. Each $\lambda^{\vee}\in X_*(T)$ determines an element in $T(\CK)$ given by $t^{\lambda^\vee}$. The assignment $\lambda^{\vee}\to G(\mathcal{O})t^{\lambda^\vee}$ is a bijection between $X_*(T)_+$ and $G(\CO)$ orbits of $G(\CK)/G(\CO)$. Denoting by $\Gr_{\lambda^\vee}$ the associated orbit, then it is a smooth quasi-projective variety since it is the quotient of an affine algebraic group by an algebraic subgroup. The reduced locus $\Gr_{G,red}$ has a stratification:
\begin{equation}
\Gr_{G,red}=\bigcup_{\lambda^\vee\in X_*(T)_+}\Gr_{\lambda^\vee}
\end{equation}
Let $\Grbar_{\lambda^\vee}$ be the Zariski closure of $\Gr_{\lambda^\vee}$. Then for $\lambda^\vee\leq \mu^\vee$ in $ X_*(T)_+$, $\Grbar_{\lambda^\vee}$ is a closed subscheme of $\Grbar_{\mu^\vee}$. This gives $\Gr_{G,red}$ an ind-scheme structure:
\begin{equation}\label{eqredstrat}
\Gr_{G,red}=\varinjlim\limits_{\lambda^\vee\in X_*(T)_+}\Grbar_{\lambda^\vee}.
\end{equation}
Each $\Grbar_{\lambda^\vee}$ is a projective variety, though usually it's very singular. In general, $\Grbar_{\lambda^\vee}$ is a normal projective variety, and it is smooth if and only if $\lambda^\vee$ is miniscule, in which case $\Grbar_{\lambda^\vee}=\Gr_{\lambda^\vee}$. These are called miniscule orbits, and are in one-to-one correspondence with the fundamental group of $G$, as well as with the number of connected components of $\Gr_G$. Note that when $G$ is semi-simple, $\Gr_G=\Gr_{G,red}$. Thus in this case, equation \eqref{eqredstrat} gives an explicit stratification of $\Gr_G$.

\subsection{The BFN space $\CBFN$}\label{BFNspace}

Now fix a finite dimensional representation $V$ of $G$. The Cartesian diagram in equation  \eqref{eqBFNdef} defines $\CBFN$ as a derived stack. In contrast to the affine Grassmannian, $\CBFN$ is not a classical ind-scheme, but a DG-indscheme. This means that it is not determined solely by its value on classical rings. Now fix an ind-scheme structure of $\Gr_G$, say $\Gr_G=\varinjlim \Gr_{G,n}$ such that each $\Gr_{G,n}$ is a projective scheme closed under the action of $\tilde{G}_\CO$. Let $G(\CK)_n$ be the pre-image of $\Gr_{G,n}$ under the projection $G(\CK)\to \Gr_G$. For each $n$, choose $N$ (that depends on $n$) large enough so that the action of $G(\CK)_n$ maps $V(\CO)$ to $z^{-N}V(\CO)$. We have the following
base-change diagram:
\begin{equation}\label{eqBFNlambda} 
\begin{tikzcd}
\mathcal{R}_{G,V,n} \rar \arrow[d] &  V(\mathcal{O})\arrow[d] \\  G(\mathcal{K})_n\times_{G(\mathcal{O})}V(\mathcal{O}) \rar &  z^{-N}V(\mathcal{O})
\end{tikzcd}
\end{equation} 
Since the bottom line of the Cartesian square in equation \eqref{eqBFNdef} is an inductive limit of the bottom line from equation \eqref{eqBFNlambda}, $\CBFN$ has the following presentation as an ind-scheme:
\begin{equation}
\CBFN=\varinjlim \mathcal{R}_{G,V,n}.
\end{equation}
Each $\mathcal{R}_{G,V,n}$ is a coconnective DG scheme, as $V(\CO)$ is a finite codimensional vector subspace in $z^{-N}V(\CO)$, and for $n\leq m$, the map $\mathcal{R}_{G,V,n}\to \mathcal{R}_{G,V,m}$ is a closed embedding. 

We will also need a local description of $\CBFN$. Let $L^-G$ be the group ind-scheme associating to an algebra $R$ the set $L^-G(R)=G(R[z^{-1}])$, and let $L^{<0}G$ be the kernel of $L^-G\to G$ sending $z^{-1}\mapsto 0$. Then according to \cite{beauville1994conformal, zhu2016introduction}, the map:
\begin{equation}
L^{<0}G\times G(\CO)\to G(\CK)
\end{equation}
is an open embedding. Thus $L^{<0}G$ is an open neighborhood of identity coset in $\Gr_G$. Over $L^{<0}G$, the vector bundle $G(\CK)\times_{G(\CO)}V(\CO)$ trivializes to $L^{<0}G\times V(\CO)$, and so over this local chart, $\CBFN$ can be represented by a pro-DG-algebra whose underlying pro-algebra is the pro-algebra of functions on the following ind-scheme:
\begin{equation}
L^{<0}G\times V(\CO)\times V(\CK)/V(\CO)[-1],
\end{equation}
and whose differential $D$ is induced from the action map:
\begin{equation}
L^{<0}G\times V(\CO) \to V(\CK)\to V(\CK)/V(\CO).
\end{equation}

\section{The Poisson Algebra $\mathcal{V}_{G,V}$}\label{3}

As explained in  \cite{oh2020poisson} and \cite{butson2021equivariant}, the algebra of local operators of the HT twist of a 4d $\mathcal{N}=2$ gauge theory has the structure of a Poisson vertex algebra, which we denote by $\mathcal{V}_{G,V}$. Here we will recall their construction. Consider the commutative Poisson vertex algebra $\mathcal{V}_{\beta\gamma-bc}$ generated by bosonic fields $(\beta,\gamma)$ with conformal weight $\frac{1}{2}$ and cohomological degree $0$, valued in the representations $V$ and $V^*$, as well as fermionic fields $(b,c)$ with conformal weight $(1,0)$ and cohomological degree $(-1,1)$, valued in the Lie algebra $\mathfrak{g}$ of $G$. The nontrivial Poisson brackets are given by:
\begin{equation}
\{\beta,\gamma\}\propto \mathrm{id}_{V},~~\{b,c\}\propto C_2(\mathfrak{g}).
\end{equation}
There is a BRST operator $Q$ defined by the current:
\begin{equation}
J_{BRST}=\mathrm{Tr}(bcc)-\beta c\gamma.
\end{equation}
The action of $Q$ is given by $Q=\{J_{BRST},-\}$ and satisfies $Q^2=0$. The Poisson algebra $\mathcal{V}_{G,V}$ is defined as the $Q$-cohomology of $\mathcal{V}_{\beta\gamma-bc}$. 

Let us now describe the vacuum module of the vertex algebra $\mathcal{V}_{G,V}$ in more detail. In fact, we will describe the vacuum module of the DG Poisson vertex algebra $(\mathcal{V}_{\beta\gamma-bc},Q)$. The vacuum module of $(\mathcal{V}_{\beta\gamma-bc},Q)$ is generated by a vacuum vector $\vert 0\rangle$ such that the positive modes acts trivially, and non-positive modes act freely. This means that, as a vector space, $\mathcal{V}_{\beta\gamma-bc}$ is given by:
\begin{equation}
\mathcal{V}_{\beta\gamma-bc}=\mathbb{C}[\beta_{k-1/2},\gamma_{k-1/2},b_{k-1},c_k]_{k\leq 0}\vert 0\rangle,
\end{equation}
with the differential $Q$ given as above. If we shift the loop weight of $V(\CO)$ by $q^{1/2}$, then the above can be identified as the following vector space:
\begin{equation}
\mathbb{C}[V(\CO)]\otimes \mathbb{C}[V^*(\CO)]\otimes \bigwedge \!\!{}^* \mathfrak{g}(\CK)/\mathfrak{g}(\CO)\otimes \bigwedge \!\!{}^*\mathfrak{g}(\CK)/z\mathfrak{g}(\CO).
\end{equation}
To understand the differential, we identify the Lie algebra of $\mathfrak{g}$ with its dual using a killing form, and view $c$ as valued in $\mathfrak{g}(\CK)^*$; then we have the following vector space:
\begin{equation}\label{eqpoissondualmatter}
\mathbb{C}[V(\CO)]\otimes \mathbb{C}[V^*(\CO)]\otimes \bigwedge \!\!{}^* \mathfrak{g}(\CK)/\mathfrak{g}(\CO)\otimes \bigwedge \!\!{}^* \left(\mathfrak{g}(\CO)\right)^*,
\end{equation}
such that the $\beta c\gamma$ part of the differential is induced by the moment map, and $\mathrm{Tr}(bcc)$ part of the differential is identified with the Chevalley-Eilenberg differential. 
The vacuum module $\mathcal{V}_{G,V}$ as a DG algebra is then identified with equation \eqref{eqpoissondualmatter} together with a differential coming from a combination of derived symplectic reduction and Chevalley-Eilenberg differential. Here $\bigwedge \!\!{}^* \left(\mathfrak{g}(\CO)\right)^*$ should be understood as the direct limit of the Chevalley-Eilenberg complex of $\mathfrak{g}(\CO)/z^m\mathfrak{g}(\CO)$. 

If we consider pure gauge theory, then the only differential is the Chevalley-Eilenberg differential, and we obtain the vector space:
\begin{equation}\label{eqpoissondualpure}
\bigwedge \!\!{}^* \mathfrak{g}(\CK)/\mathfrak{g}(\CO)\otimes \bigwedge \!\!{}^* \left(\mathfrak{g}(\CO)\right)^*
\end{equation}
with the CE differential. 

Physically, the space of local operators is not simply the cohomology of $(\mathcal{V}_{\beta\gamma-bc},Q)$. This is due to the fact that when computing correlation functions, the ghosts should not appear as initial or final states. Mathematically, this amounts to, after taking cohomology, projecting to  $c_0$-ghost-number zero. This is the same as taking invariants with respect to the Lie group $G$ (the constant gauge transformations) by hand, instead of derived invariants of its Lie algebra. We denote the resulting space by $\pi_0\mathcal{V}_{G,V}$. From the above discussions, we have a quasi-isomorphism:
\begin{equation}\label{eqpoissonmatter}
\pi_0\mathcal{V}_{G,V}\cong \left[\mathbb{C}[V(\CO)]\otimes \mathbb{C}[V^*(\CO)]\otimes \bigwedge \!\!{}^* \mathfrak{g}(\CK)/\mathfrak{g}(\CO)\otimes  \bigwedge \!\!{}^* (z\mathfrak{g}(\CO))^*\right]^G.
\end{equation}
Note that $\mathfrak{g}(\CK)/\mathfrak{g}(\CO)$ and $(z\mathfrak{g}(\CO))^*$ contribute the same factor to the Euler character. However, their roles are not symmetric, since one of them is used for symplectic reduction and the other is for derived group invariants. This difference will show up in the geometric computation as well. 

\vspace{10pt}
\noindent\textbf{Remark.} The Poisson vertex algebra $\mathcal{V}_{G,V}$ exists for any gauge theory. However, when the gauge theory is super-conformal, which happens when $C_2(V)=C_2(G)$, this algebra has a deformation through the work of \cite{oh2019chiral} and  \cite{butson2021equivariant}, and the deformed algebra is identified with the conformal vertex algebra (VOA) first studied in \cite{beem2015infinite}. Their construction is as follows: the algebra $\mathcal{V}_{\beta\gamma-bc}$ has a deformation quantization into the VOA $V^{\hbar}_{\beta\gamma-bc}$ generated by bosonic fields $\beta,\gamma$ and fermionic fields $b,c$ with OPE:
\begin{equation}
\gamma(z)\beta(w)\sim \frac{\hbar \mathrm{id}_V}{z-w},~b(z)c(w)\sim \frac{\hbar C_2(G)}{z-w}.
\end{equation}
The action of $Q$ is promoted to the action of $Q_{BRST}$ via:
\begin{equation}
Q_{BRST}\CO(z)=\oint_w J_{BRST}(z+w)\CO(z). 
\end{equation}
It squares to zero precisely when $C_2(V)=C_2(G)$, and the cohomology of $Q_{BRST}$ gives the deformation quantization of $\mathcal{V}_{G,V}$. In this case, it is expected that in the category of line operators, the Schur functor is trivial, or in other words, the left dual of a line operator is isomorphic to its right dual. We will not prove this in our paper. This is one direct consequence of superconformal symmetry in the category of line operators.

\section{Geometric Computations of the Space of Local Operators}\label{4}

Let $G$ be a reductive Lie group and $V$ a finite-dimensional representation of $G$. In this section, we will specify the category of line operators in $T_{G,V}$ using $\infty$-category language, identify the identity line operator, and compute its endomorphism algebra. This will be done in steps. In Section \ref{Categoryofsheaves}, we will introduce the category of line operators using the machinery of \cite{raskin2020homological}. In Section \ref{Homspaces} we will briefly introduce the Hom functor between two line operators. In Section \ref{PuregaugessG}, we obtain the endomorphism algebra for pure gauge theory($V=0$) in Theorem \ref{Extidentitypure}. In Section \ref{AbelianGaugeGroup}, as a simple yet nontrivial example, we look at pure abelian gauge theory.  In Section \ref{gaugetheorywithmatter}, we obtain the endomorphism algebra for gauge theories with matter in Theorem \ref{Extidentitymatter}. In both cases, we find agreement with the Poisson vertex algebra $\pi_0 \mathcal{V}_{G,V}$ introduced in Section \ref{3}. We comment that in this section,  we will write $\mathrm{Sym}^{\bullet}\!(V)$ for the free graded-commutative algebra generated by a graded vector space $V$. In particular, when $V$ is in even degree, then it is a polynomial algebra; when $V$ is in odd degree, then it is an exterior algebra. 

\subsection{The Category of Equivariant Coherent Sheaves}\label{Categoryofsheaves}

A reasonable DG ind-scheme, as defined in \cite[Definition 6.8.1]{raskin2020homological}, is a convergent prestack $X$ such that $X=\varinjlim X_i$ such that each $X_i$ is quasi-compact, quasi-separated and eventually coconnective, and that $X_i\to X_{j}$ is  almost finitely-presented closed embeddings. $\Gr_G$ and $\CBFN$ are examples of such reasonable DG ind-schemes.  Let $H$ be a classical affine group scheme that acts on $X$. Then the quotient stack $X/H$ is called a weakly renormalizable pre-stack following \cite[Definition 6.28.1]{raskin2020homological}, and one can define the category $\mathrm{IndCoh}^*(X/H)$ via a right Kan extension:
\begin{equation}
\IndCoh^*(X/H):= \lim\limits_{f: S\to X/H \text{ flat}} \IndCoh^*(S),
\end{equation}
where the limit is taken over all reasonable DG ind-schemes flat over $X/H$, using the functoriality of $f^{*,\IndCoh}$. We have the following equivalence:
\begin{equation}
\IndCoh^*(X/H)\cong \IndCoh^*(X)^{H,w,naive}:=\mathrm{Hom}_{H-\mathrm{mod}_{weak,naive}}(\mathrm{Vect}, \IndCoh^*(X)),
\end{equation}
where the right hand side is the naive weakly equivariant category with respect to the action of $H$ as defined in \cite{raskin2020homological} Section 5. This in particular, may not be equivalent to the ind-completion of its compact object. This category may seem abstract, but one can unpack it using flat descent. Recall that given a flat cover $T\to S$, one can consider the associated Cech nerve:
\begin{equation}
T^{\times_{S}^{*+1}}.
\end{equation}
Applying this to the flat cover $X\to X/H$, the Cech nerve is:
\begin{equation}\label{eqcosimplicial}
X^{\times_{X/H}^{*+1}}= X\lstack{3} X\times H \lstack{5} X\times H \times H\cdots
\end{equation}
By \cite[Theorem 6.25.1]{raskin2020homological}:
\begin{equation}\label{eqcosimplicialcat}
\IndCoh^*(X/H)\cong \mathrm{Tot}_{semi}(\IndCoh(X^{\times_{X/H}^{*+1}})).
\end{equation}
The right hand side is a semi-simplicial set of categories that only involve categories of sheaves on ind-schemes. Suppose further that $H$ acts on each $X_i$ such that $X/H=\varinjlim X_i/H$, we can write:
\begin{equation}
\mathrm{Tot}_{semi}(\IndCoh(X^{\times_{X/H}^{*+1}}))=\mathrm{Tot}_{semi}(\lim\limits_{\text{upper-}!}\IndCoh(X_i\times_{X/H}X^{\times_{X/H}^{*+1}})).
\end{equation}
Commuting the limit on the right hand side using \cite[Lemma 6.17.2]{raskin2020homological}, noticing that $X_i\times_{X/H}X^{\times_{X/H}^{*+1}})$ is the Cech nerve of $X_i\to X_i/H$, we get:
\begin{equation}
\mathrm{Tot}_{semi}(\IndCoh(X^{\times_{X/H}^{*+1}}))=\lim\limits_{\text{upper-}!} \IndCoh(X_i/H).
\end{equation}
By \cite[Lemma 1.3.3]{gaitsgory2012notes}, we may change the limit over upper-$!$ to the colimit over lower-$*$:
\begin{equation}\label{eqinjquotientindcoh}
\IndCoh^*(X/H)\cong \varinjlim_{\text{lower-}*} \IndCoh^*(X_i/H).
\end{equation}

Now we specialize this story to the BFN space. Let $\mathbb{C}^*$ act as the two-fold cover of the loop rotation. Both $\Gr_G$ and $\CBFN$ have an action of $G(\CO)\rtimes \mathbb{C}^*$.\footnote{We will use the two fold cover of the group of loop rotations here, since this will allow us to shift gradings by $q^{1/2}$, which is necessary for matching with the physical Schur indices.} We will denote this group by $\tilde{G}_\CO$. From the above discussion, we can define categories $\IndCoh(\tilde{G}_\CO\setminus \Gr_G)$ and more generally, $\IndCoh(\tilde{G}_\CO\setminus \CBFN)$. Moreover, if we fix a stratification $\{\mathcal{R}_{G,V, n}\}$ of $\CBFN$ as in Section \ref{BFNspace}, then:
\begin{equation}
\IndCoh(\tilde{G}_\CO\setminus \CBFN)\cong \varinjlim \IndCoh(\tilde{G}_\CO\setminus \mathcal{R}_{G,V,n}).
\end{equation}
For each $n$, the category of coherent sheaves $\Coh(\tilde{G}_\CO\setminus \mathcal{R}_{G,V,n})$ is the full subcategory of $\IndCoh(\tilde{G}_\CO\setminus \mathcal{R}_{G,V,n})$ consisting of objects whose pull-back to $\mathcal{R}_{G,V,n}$ is coherent. The category of equivariant coherent sheaves on $\CBFN$, $\Coh(\tilde{G}_\CO\setminus \CBFN)$, is defined as the full subcategory of $\IndCoh(\tilde{G}_\CO\setminus \CBFN)$ whose objects are the images of $\Coh(\tilde{G}_\CO\setminus \mathcal{R}_{G,V,n})$ under the above colimit. This category $\Coh(\tilde{G}_\CO\setminus \CBFN)$ is expected to be the category of line operators for the theory $T_{G,V}$, and the derived Hom between objects in this category is expected to be the space of local operators at the junction of two lines. 

\subsection{On Hom Spaces}\label{Homspaces}

To obtain the space of local operators at the junction of two line operators, we need to take the Hom space between two line operators as a DG vector space. Intuitively, this is the derived Hom spaces between two coherent sheaves. In this section, we will briefly introduce the setting in which this enriched Hom can be taken. Let ${\mathcal C}$ be a presentable monoidal DG category and ${\mathcal M}$ be a presentable DG module category of~${\mathcal C}$. For any pair of objects $(M_1,M_2)$ in ${\mathcal M}$,  the object
\begin{equation}
\mathrm{Hom}^{\mathcal C}(M_1,M_2)\in {\mathcal C}
\end{equation}
is defined by the following adjunction property:
\begin{equation}
\mathrm{Hom}_{\mathcal M}(-\otimes M_1, M_2)=\mathrm{Hom}_{\mathcal C}(-, \mathrm{Hom}^{\mathcal C}(M_1,M_2)).
\end{equation}
Let $X$ be a reasonable DG indscheme acted on by a smooth affine group scheme $H$. Then $\mathrm{IndCoh}(X/H)$ is a module category over the monoidal category $\mathrm{IndCoh}(\mathbb{B}H)\cong \mathrm{QCoh}(\mathbb{B}H)$, where $\mathbb{B}H=\mathrm{pt}/H$ is the classifying stack of $H$. We thus obtain a Hom functor:
\begin{equation}
\mathrm{Hom}^{\mathrm{QCoh}(\mathbb{B}H)}\!(-,-): \mathrm{IndCoh}(X/H)^{op}\times \mathrm{IndCoh}(X/H)\to \mathrm{QCoh}(\mathbb{B}H). 
\end{equation}
We will abbreviate this by $\mathrm{Hom}^{\mathbb{B}H}$. Specify this to our setting, we have the Hom functor:
\begin{equation}
\begin{aligned}
\mathrm{Hom}^{\mathbb{B}\mathbb{C}^*}\!\!\!\!: ~~~ \IndCoh^*(\tilde{G}_\CO\setminus \CBFN)^{op}&\times \IndCoh^*(\tilde{G}_\CO\setminus \CBFN) \\ & \longrightarrow \mathrm{QCoh}(\mathbb{B}\mathbb{C}^*)\cong \mathrm{IndCoh}(\mathbb{B}\mathbb{C}^*).
\end{aligned}
\end{equation}
This will be the main player of this section for many of the computations. We will write $\mathrm{End}^{\mathbb{B}\mathbb{C}^*}$ if the two arguments of $\mathrm{Hom}$ are identical. 

\vspace{10pt}
\noindent\textbf{Remark.} This definition of $\mathrm{Hom}$ spaces seem to be abstract, but in our example it is a concrete one: first of all, it will be given by a colimit of $\mathrm{Hom}$ spaces computed on each closed orbit of $G(\CO)$; secondly, on each orbit, it is the usual dg vector space of $\mathrm{Hom}$, which can be computed by choosing an injective resolution of the second argument. 

\subsection{Bulk Local Operators in Pure Gauge Theory}\label{PuregaugessG}

In pure-gauge, the category of line operators is $\mathrm{Coh}(\tilde{G}_\CO\setminus\Gr_G)$. This category is a monoidal category, with monoidal unit given by $\CO_{[e]/ \tilde{G}_\CO}$, the structure sheaf of the identity coset $[e]$ with the trivial $\tilde{G}_\CO$ equivariant structure. Our goal is to compute the space:
\begin{equation}
\mathrm{End}^{\mathbb{B}\mathbb{C}^*}(\CO_{[e]/ \tilde{G}_\CO})
\end{equation}
as a $\mathbb{C}^*$-DG vector space. The remainder of this section is devoted to the computation of this space, up to quasi-isomorphism. The idea of the computation is the following:
\begin{itemize}
\item First, prove that one can factor the computation into two steps: computing $\mathrm{End}^{\mathbb{B}\tilde{G}_\CO}(\CO_{[e]/ \tilde{G}_\CO})$; then taking the (derived-)invariant subspace with respect to the $G(\CO)$ action.

\item Computing the derived $G(\CO)$ invariants using the Chevalley-Eilenberg cochain complex.

\item Computing  $\mathrm{End}^{\mathbb{B}\tilde{G}_\CO}(\CO_{[e]/ \tilde{G}_\CO})$ using formal completion. 

\end{itemize}

The result, stated in Theorem \ref{Extidentitypure}, coincides with $\pi_0 \mathcal{V}_{G,0}$ after some appropriate degree shift. 

\subsubsection{Decomposing the Hom Functor}

Let $H$ be a smooth affine group scheme that can be written as:
\begin{equation}
H= H_0\rtimes T
\end{equation}  
for two smooth affine group schemes $H_0$ and $T$. Assume also that $T$ is of finite type. Let $Y_n$ be finite-type classical $H$-schemes such that $Y_n\to Y_{n+1}$ are closed embeddings of $H$-schemes. Denote by $Y=\varinjlim Y_n$ and $\mathcal{Y}=\varinjlim Y_n/H$. Let $X$ be a finite-type classical $H$-scheme together with a closed-embedding of $H$-schemes $i: X\to Y$. Denote by $\mathcal{X}=X/H$. Let $(\mathcal{F},\mathcal{G})$ be a pair of objects in $\mathrm{Coh}(\mathcal{X})$. We would like to understand
\begin{equation}
\mathrm{Hom}^{\mathrm{pt}/T}(i_{*,\IndCoh}\mathcal{F},i_{*,\IndCoh}\mathcal{G}).
\end{equation} 
Denote by $\pi$ the natural projection $\mathcal{Y}\to \mathrm{pt}/H=\mathbb{B}H$, the classifying stack of $H$, and by $\pi_0$ the natural map $\mathbb{B}H\to \mathbb{B}T$. Since $\IndCoh(\mathcal{Y})$ is a module category of $\IndCoh(\mathbb{B}H)$, we have an object:
\begin{equation}
\mathrm{Hom}^{\mathbb{B}H}(i_{*,\IndCoh}\mathcal{F},i_{*,\IndCoh}\mathcal{G})\in \IndCoh(\mathrm{pt}/H).
\end{equation}
Now if we view $\IndCoh(\mathbb{B}H)$ as a module category of $\IndCoh(\mathbb{B}T)$ via the functor $\pi_0^*$, we will have an object:
\begin{equation}
\mathrm{Hom}^{\mathbb{B}T}\left(\CO_{\mathbb{B}H}, \mathrm{Hom}^{\mathbb{B}H}(i_{*,\IndCoh}\mathcal{F},i_{*,\IndCoh}\mathcal{G})\right).
\end{equation}

\newtheorem{LemHequiv}{Lemma}[section]

\begin{LemHequiv}\label{LemHequiv1}
The following is a quasi-isomorphism of $T$ modules:
\begin{equation}\label{eqcolimT}
\mathrm{Hom}^{\mathbb{B}T}(i_{*,\IndCoh}\mathcal{F},i_{*,\IndCoh}\mathcal{G})\cong \mathrm{Hom}^{\mathbb{B}T}\left(\CO_{\mathbb{B}H}, \mathrm{Hom}^{\mathbb{B}H}(i_{*,\IndCoh}\mathcal{F},i_{*,\IndCoh}\mathcal{G})\right).
\end{equation}

\end{LemHequiv}

\begin{proof}
Let $V$ be an object of  $\IndCoh(\mathbb{B}T)$, then:
\begin{equation}
\begin{aligned}
\mathrm{Hom}_{\IndCoh(\mathbb{B}T)}& \left( V,  \mathrm{Hom}^{\mathbb{B}T}\left(\CO_{\mathbb{B}H}, \mathrm{Hom}^{\mathbb{B}H}(i_{*,\IndCoh}\mathcal{F},i_{*,\IndCoh}\mathcal{G})\right) \right)\\ &\cong \mathrm{Hom}_{\IndCoh(\mathbb{B}H)}\left(\pi_0^*V,  \mathrm{Hom}^{\mathbb{B}H}(i_{*,\IndCoh}\mathcal{F},i_{*,\IndCoh}\mathcal{G})\right)\\
&\cong \mathrm{Hom}_{\IndCoh(\mathcal Y)}\left(V\otimes i_{*,\IndCoh}\mathcal{F}, i_{*,\IndCoh}\mathcal{G}\right)\\ & \cong \mathrm{Hom}_{\IndCoh(\mathbb{B}T)}\left(V, \mathrm{Hom}^{\mathbb{B}T}\left(i_{*,\IndCoh}\mathcal{F},i_{*,\IndCoh}\mathcal{G}\right)\right). 
\end{aligned}
\end{equation} 
This proves the claim.

\end{proof}

This statement says that we can first compute the endomorphism of $i_{*,\IndCoh}\mathcal{F}$ and $i_{*,\IndCoh}\mathcal{G}$ as an $H$-module, and then compute invariants with respect to $H_0$. However, this is not the best way to understand this Hom space, since $\IndCoh(\mathbb{B}H)$ is not compactly generated. In \cite[Section 5.11]{raskin2020homological}, the author defined another category that is compactly generated. Denote by $\mathrm{Rep}(H)^c$ the monoidal subcategory of $\IndCoh(\mathbb{B}H)$ consisting of objects whose images under the forgetful functor $\IndCoh(\mathbb{B}H)\to \mathrm{Vect}$ are compact, and $\mathrm{Rep}(H)=\mathrm{Ind}(\mathrm{Rep}(H)^c)$, the ind-completion. This category is compactly generated, and if $H$ is a smooth affine algebraic group, then it is equivalent to $\IndCoh(H)$. In particular, $\IndCoh(\mathbb{B}T)\cong \mathrm{Rep}(T)$. 

 Moreover, by \cite[Lemma 5.16.2]{raskin2020homological}, if $H=\lim H_i$ for $H_i$ finite dimensional smooth algebraic groups, then $\mathrm{Rep}(H)=\varinjlim \mathrm{Rep}(H_i)$, and so the understanding of $\mathrm{Rep}(H)$ can be reduced to understanding representations of finite-dimensional algebraic groups. 
 
Since the action of $\mathrm{Rep}(T)$ on both $\IndCoh(\mathbb{B}H)$ and $\IndCoh(\mathcal Y)$ factors through an action of $\mathrm{Rep}(H)$, we can modify Lemma \ref{LemHequiv1} into the following:
\newtheorem{LemHRep}[LemHequiv]{Lemma}

\begin{LemHRep}\label{LemHRep}
The following is a quasi-isomorphism of $T$ modules:
\begin{equation}\label{eqcolimT}
\mathrm{Hom}^{\mathbb{B}T}(i_{*,\IndCoh}\mathcal{F},i_{*,\IndCoh}\mathcal{G})\cong \mathrm{Hom}^{\mathbb{B}T}\left(\CO_{\mathbb{B}H}, \mathrm{Hom}^{\mathrm{Rep}(H)}(i_{*,\IndCoh}\mathcal{F},i_{*,\IndCoh}\mathcal{G})\right).
\end{equation}

\end{LemHRep}

The object $\mathrm{Hom}^{\mathrm{Rep}(H)}(i_{*,\IndCoh}\mathcal{F},i_{*,\IndCoh}\mathcal{G})$ behaves better with colimit in equation \eqref{eqinjquotientindcoh}:
\newtheorem{colimRepH}[LemHequiv]{Proposition}

\begin{colimRepH}\label{colimRepH}
Denote by $\mathcal{F}_k$ and $\mathcal{G}_k$ the pushforward of $\mathcal{F}$ and $\mathcal G$ to $Y_k/H$. There is a qausi-isomorphism in $\mathrm{Rep}(H)$:
\begin{equation}\label{eqcolimTreal}
\mathrm{Hom}^{\mathrm{Rep}(H)}(i_{*,\IndCoh}\mathcal{F},i_{*,\IndCoh}\mathcal{G})\cong\varinjlim \mathrm{Hom}^{\mathrm{Rep}(H)}(\mathcal{F}_k,\mathcal{G}_k). 
\end{equation}
\end{colimRepH}

\begin{proof}
Given $V\in \mathrm{Rep}(H)^c$, we have:
\begin{equation}
\begin{aligned}
\mathrm{Hom}_{\mathrm{Rep}(H)}&\left(V, \mathrm{Hom}^{\mathrm{Rep}(H)}(i_{*,\IndCoh}\mathcal{F}, i_{*,\IndCoh}\mathcal{G})\right) \\ & \cong \mathrm{Hom}_{\IndCoh(\mathcal{Y})}(V\otimes i_{*,\IndCoh}\mathcal{F}, i_{*,\IndCoh}\mathcal{G})\\(\text{by equation } \eqref{eqinjquotientindcoh}) &\cong \varinjlim  \mathrm{Hom}_{\IndCoh(Y_k/H)}(V\otimes \mathcal{F}_k,\mathcal{G}_k)\\ &\cong \varinjlim \mathrm{Hom}_{\mathrm{Rep}(H)}\left(V,  \mathrm{Hom}^{\mathrm{Rep}(H)}(\mathcal{F}_k,\mathcal{G}_k)\right) \\(\text{since } V \text{ is compact}) &\cong \mathrm{Hom}_{\mathrm{Rep}(H)}\left(V, \varinjlim \mathrm{Hom}^{\mathrm{Rep}(H)}(\mathcal{F}_k,\mathcal{G}_k)\right)
\end{aligned}
\end{equation}
Since $\mathrm{Rep}(H)$ is compactly generated, this proves the claim. 
\end{proof}
 
 The object $\mathrm{Hom}^{\mathrm{Rep}(H)}(\mathcal{F}_k,\mathcal{G}_k)$ may seem to be abstract at first, but we can show that this is a familiar object: the underlying vector space of this object is the derived Hom between $\mathcal{F}_k$ and $\mathcal{G}_k$ as sheaves over $Y_k$. Denote by $\mathbf{Oblv}$ the forgetful functor $\mathrm{Rep}(H)\to \mathrm{Vect}$. This is the composition of $\Psi: \mathrm{Rep}(H)\to \IndCoh(\mathbb{B}H)$ with the forgetful functor $ \IndCoh(\mathbb{B}H)\to \mathrm{Vect}$. Denote also by $p_k$  the projection $X_k\to X_k/H$. We claim:
 \newtheorem{LemHequiv2}[LemHequiv]{Proposition}
 
 \begin{LemHequiv2}\label{LemHequiv2}
There is a quasi-isomorphism: 
\begin{equation}
\mathbf{Oblv}\mathrm{Hom}^{\mathrm{Rep}(H)}(\mathcal{F}_k,\mathcal{G}_k)\cong \mathrm{Hom}^{\mathrm{Vect}}(p_k^*\mathcal{F}_k,p_k^*\mathcal{G}_k).
\end{equation}
Here the left hand side of the above equation is the underlying DG vector space of $\mathrm{Hom}^{\mathrm{Rep}(H)}(\mathcal{F}_k,\mathcal{G}_k)$.
\end{LemHequiv2}

We need the following Lemma:

\newtheorem{Lemma1}[LemHequiv]{Lemma}

\begin{Lemma1}\label{Lemma1}
Let $p_0: \mathrm{pt}\to \mathbb{B}H$ be the projection, then:
\begin{equation}
p_0^*\mathrm{Hom}^{\mathbb{B}H}(\mathcal{F}_k,\mathcal{G}_k)\cong  \mathrm{Hom}^{\mathrm{Vect}}(p_k^*\mathcal{F}_k,p_k^*\mathcal{G}_k).
\end{equation}
\end{Lemma1}

 \begin{proof}
Denote by $\pi_k$ the projection $Y_k/H\to \mathbb{B}H$. We know that $\IndCoh(Y_k/H)$ is a module category of $\mathrm{QCoh}(Y_k/H)$, and this is compatible with the monoidal functor:
\begin{equation}
\pi_k^*: \IndCoh(\mathbb{B}H)\cong \mathrm{QCoh}(\mathbb{B}H)\to \mathrm{QCoh}(Y_k/H).
\end{equation} 
Using adjunction property, we have:
\begin{equation}
\mathrm{Hom}^{\mathbb{B}H}(\mathcal{F}_k,\mathcal{G}_k)=\mathrm{Hom}^{\mathbb{B}H}\left(\mathcal{O}_{Y_k/H},\mathrm{Hom}^{\mathrm{QCoh}(Y_k/H)}(\mathcal{F}_k,\mathcal{G}_k)\right).
\end{equation}
The right hand side of the above equation can be identified with:
\begin{equation}
(\pi_k)_* \mathrm{Hom}^{\mathrm{QCoh}(Y_k/H)}(\mathcal{F}_k,\mathcal{G}_k).
\end{equation}
We are thus interested in $p_0^*(\pi_k)_* \mathrm{Hom}^{\mathrm{QCoh}(X_k/H)}(\mathcal{F}_k,\mathcal{G}_k)$. Consider now the Cartesian diagram:
\begin{equation}
\begin{tikzcd}
X_k \rar{\tilde{\pi}_k} \arrow[d, "p_k"] &  \mathrm{pt} \arrow[d, "p_0"] \\ X_k/H \rar{\pi_k} &  \mathrm{pt}/H
\end{tikzcd}
\end{equation} 
Using base-change property of $\mathrm{QCoh}$, we obtain:
\begin{equation}\label{eqLemHequiv2-1}
p_0^*(\pi_k)_* \mathrm{Hom}^{\mathrm{QCoh}(X_k/H)}(\mathcal{F}_k,\mathcal{G}_k)\cong (\tilde{\pi}_k)_* p_k^*\mathrm{Hom}^{\mathrm{QCoh}(X_k/H)}(\mathcal{F}_k,\mathcal{G}_k). 
\end{equation}
Now by \cite[Proposition 9.5.3.3]{lurie2018spectral}:
\begin{equation}
p_k^*\mathrm{Hom}^{\mathrm{QCoh}(X_k/H)}(\mathcal{F}_k,\mathcal{G}_k)\cong \mathrm{Hom}^{\mathrm{QCoh}(X_k)}(p_k^*\mathcal{F}_k,p_k^*\mathcal{G}_k).
\end{equation}
Putting this into equation \eqref{eqLemHequiv2-1} we obtain the desired result.

 \end{proof}
 
 \begin{proof}[Proof of Proposition \ref{LemHequiv2}.] 
 By adjunction property, there is a quasi-isomorphism:
 \begin{equation}
 \mathrm{Hom}^{\mathrm{Rep}(H)}(\mathcal{F}_k,\mathcal{G}_k)\cong \mathrm{Hom}^{\mathrm{Rep}(H)}\left(\CO_{\mathbb{B}H}, \mathrm{Hom}^{\mathbb{B}H}(\mathcal{F}_k,\mathcal{G}_k)\right).
 \end{equation}
 Since $Y_k$ is a classical finite-type scheme and $\mathcal{F}_k$ and $\mathcal{G}_k$ are coherent, by Lemma \ref{Lemma1},  the Hom space $\mathrm{Hom}^{\mathbb{B}H}(\mathcal{F}_k,\mathcal{G}_k)$ is an object in $\IndCoh(\mathbb{B}H)^+$, which is equivalent to $\mathrm{Rep}(H)^+$ via $\Psi$. Thus:
 \begin{equation}
 \Psi  \mathrm{Hom}^{\mathrm{Rep}(H)}(\mathcal{F}_k,\mathcal{G}_k)\cong \mathrm{Hom}^{\mathbb{B}H}(\mathcal{F}_k,\mathcal{G}_k).
 \end{equation}
 Since $\mathbf{Oblv}=p_0^*\circ \Psi$, this and Lemma \ref{Lemma1} gives the desired result. 
 \end{proof}
 
\noindent\textbf{Remark.} The above discussions suggest that the sheaf $\mathrm{Hom}^{\mathrm{QCoh}(Y_k/H)}(\mathcal{F}_k,\mathcal{G}_k)$ is the usual Hom sheaf between $\mathcal{F}_k$ and $\mathcal{G}_k$ on $Y_k$ with the canonical $H$ equivariant structure. The( derived) global section of this sheaf over $Y_k$ as an $H$ module is identified with $\mathrm{Hom}^{\mathrm{Rep}(H)}(\mathcal{F}_k,\mathcal{G}_k)$. The $H$ module $\mathrm{Hom}^{\mathrm{Rep}(H)}(i_{*,\IndCoh}\mathcal{F},i_{*,\IndCoh}\mathcal{G})$ is the colimit of $\mathrm{Hom}^{\mathrm{Rep}(H)}(\mathcal{F}_k,\mathcal{G}_k)$. 

\vspace{10pt}

We will apply this to the affine Grassmannian $\Gr_G$. Fix a stratification $\Gr_{G}=\varinjlim \Gr_{G,n}$ such that $\Gr_{G,n}$ is a projective scheme closed under the action of $\tilde{G}_\CO$. Take $\mathcal{F}$ and $\mathcal{G}$ to be objects in $\Coh(\tilde{G}_\CO\setminus \Gr_{G,n})$, viewed as objects in $\Coh(\tilde{G}_\CO\setminus \Gr_G)$. Lemma \ref{LemHRep} implies:
\begin{equation}
\mathrm{Hom}^{\mathbb{B}\mathbb{C}^*}(\mathcal{F},\mathcal{G})\cong \mathrm{Hom}^{\mathbb{B}\mathbb{C}^*}\left(\CO_{\mathbb{B}\tilde{G}_\CO}, \mathrm{Hom}^{\mathrm{Rep}(\tilde{G}_\CO)}(\mathcal{F},\mathcal{G})\right).
\end{equation}
Proposition \ref{colimRepH} shows that the $\tilde{G}_\CO$ module $\mathrm{Hom}^{\mathrm{Rep}(\tilde{G}_\CO)}(\mathcal{F},\mathcal{G})$ is a colimit of Hom spaces on finite-dimensional strata $\Gr_{G,k}$, namely $\mathrm{Hom}^{\mathrm{Rep}(\tilde{G}_\CO)}(\mathcal{F}_k,\mathcal{G}_k)$. These are bounded from below independent of $k$ by Proposition \ref{LemHequiv2}, and so can be identified with $\mathrm{Hom}^{\mathbb{B}H}(\mathcal{F}_k,\mathcal{G}_k)$. The functor $\mathrm{Hom}^{\mathbb{B}\mathbb{C}^*}(\CO_{\mathbb{B}\tilde{G}_\CO},-)$ thus computes the derived invariants of these $\tilde{G}_\CO$ modules with respect to the normal subgroup $G(\CO)$. Let $G_{>n}$ be the normal subgroup defined by $G(1+z^n\CO)$, then $\tilde{G}_\CO\cong \lim \tilde{G}_\CO/G_{>n}$, and so by \cite[Lemma 5.16.2]{raskin2020homological}, $\mathrm{Rep}(\tilde{G}_\CO)\cong \varinjlim \mathrm{Rep}(\tilde{G}_\CO/G_{>n})$, and so taking $G(\CO)$ invariants of modules in $\mathrm{Rep}(\tilde{G}_\CO)$ can be calculated by analyzing invariants of finite algebraic groups. This is what we turn to next. 

\subsubsection{Equivariance with Respect to $G(\CO)$}

In this section, we will deal with $G(\CO)$ invariants. Let $G_{>0}$ be the kernel of the group homomorphism $G(\CO)\to G$ given by mapping $g[z]\to g(0)$, and let $\mathfrak{g}_{>0}$ be it's Lie algebra. For each representation $V$, there is an associated Chevalley-Eilenberg cochain complex:
\begin{equation}
V\otimes   \mathrm{Sym}^{\bullet}\!\!\left(\left(z\mathfrak{g}\left(\CO\right)\right)^*[-1]\right),
\end{equation}
in which the differential $V\to V\otimes (z\mathfrak{g}(\CO))^*$ is induced by the action of $z\mathfrak{g}(\CO)$ on $V$. In this section, we will prove:

\newtheorem{PropGOequiv}[LemHequiv]{Proposition}

\begin{PropGOequiv}\label{PropGOequiv}
 Let $V$ be an algebraic $\tilde{G}_\CO$ representation, there is a quasi-isomorphism:
\begin{equation}\label{eqGOequiv}
\mathrm{Hom}^{\mathbb{B}\mathbb{C}^*}(\mathcal{O}_{\mathbb{B}\tilde{G}_\CO}, V)\cong \left[ V\otimes \mathrm{Sym}^{\bullet}\!\!\left(\left(z\mathfrak{g}\left(\CO\right)\right)^*[-1]\right)\right]^G.
\end{equation}
Here $[-]^G$ means taking the $G$ invariant part of a representation. 

\end{PropGOequiv}

\begin{proof}
We have a short exact sequence of groups:
\begin{equation}
1\to G_{>0}\to G(\CO)\to G\to 1,
\end{equation}
which gives a natural equivalence of functors:
\begin{equation}
\mathrm{Hom}^{\mathbb{B} \mathbb{C}^*}\!\!(\mathcal{O}_{\mathbb{B}\tilde{G}_\CO}, V)\cong \mathrm{Hom}^{\mathbb{B}\mathbb{C}^*}\!\!\left(\mathcal{O}_{\mathbb{B}(G \times \mathbb{C}^*)}, \mathrm{Hom}^{\mathbb{B}(G\times \mathbb{C}^*)}(\mathcal{O}_{\mathbb{B}\tilde{G}_\CO}, V)\right)
\end{equation}
Since $G$ is reductive, the category of algebraic representations of $G$ is semi-simple, which implies that:
\begin{equation}
 \mathrm{Hom}^{\mathbb{B}\mathbb{C}^*}\!\!\left(\mathcal{O}_{\mathbb{B}(G \times \mathbb{C}^*)}, \mathrm{Hom}^{\mathbb{B}(G\times \mathbb{C}^*)}(\mathcal{O}_{\mathbb{B}\tilde{G}_\CO}, V)\right)=\left[\mathrm{Hom}^{\mathbb{B}(G\times \mathbb{C}^*)}(\mathcal{O}_{\mathbb{B}\tilde{G}_\CO}, V)\right]^G,
\end{equation}
where $[-]^G$ is taking ordinary $G$ invariants. Thus we only need to understand $\mathrm{Hom}^{\mathbb{B}(G\times \mathbb{C}^*)}(\mathcal{O}_{\mathbb{B}\tilde{G}_\CO}, V)$, which is the( derived) $G_{>0}$ invariants of $V$. To understand this, we need:
\newtheorem{LemGnilp}[LemHequiv]{Lemma}

\begin{LemGnilp}\label{LemGnilp}
Let $K$ be a finite-dimensional simply-connected unipotent Lie group and $V$ be an algebraic representation of $K$, then $\mathrm{RHom}(\mathbb{C}, V)\cong \mathrm{H}^*\left(V\otimes \mathrm{Sym}^{\bullet}\!\!\left(\mathfrak{k}^*[-1]\right)\right)=: \mathcal{H}^*( \mathfrak{k},V)$, where $V\otimes \mathrm{Sym}^{\bullet}\!\!\left(\mathfrak{k}^*[-1]\right)$ is the Chevalley Eilenberg cochain complex of $V$ as a $\mathfrak{k}$ module. 

\end{LemGnilp}

Let us assume this for now and apply it to $G_{>0}$. From $\mathrm{Rep}(G_{>0})=\varinjlim\mathrm{Rep}(G_{>0}/G_{>m})$, we see that for any finite-dimensional representation $V$:
 \begin{equation}
 \mathrm{RHom}_{G_{>0}}(\mathbb{C},V)=\varinjlim\limits_m\mathrm{RHom}_{G_{>0}/G_{>m}}(\mathbb{C},V).
 \end{equation} 
 Now the Lie group $G_{>0}/G_{>m}$ is unipotent simply-connected, whose Lie algebra is $z\mathfrak{g}(\CO)/z^{m+1}\mathfrak{g}(\CO)$, so for finite-dimensional $V$, one has:
\begin{equation}
\mathrm{RHom}_{G_{>0}/G_{>m}}(\mathbb{C},V)\cong V\otimes \mathrm{Sym}^{\bullet}\!\!\left(\left(z\mathfrak{g}(\CO)/z^{m+1}\mathfrak{g}(\CO)\right)^*[-1]\right)
\end{equation}
Taking co-limit over $m$, one obtain, for any finite-dimensional(and more generally algebraic) representation $V$:
\begin{equation}
\mathrm{RHom}_{G_{>0}}(\mathbb{C},V)\cong V\otimes \mathrm{Sym}^{\bullet}\!\!\left(\left(z\mathfrak{g}\left(\CO\right)\right)^*[-1]\right).
\end{equation}
This completes the proof. 

\end{proof}

For completeness, we present the proof of Lemma \ref{LemGnilp} here:

\begin{proof}[Proof of Lemma \ref{LemGnilp}]
Clearly $\mathcal{H}^0=\mathrm{Hom}_K(\mathbb{C},V)$, so by the usual idea of homological algebra(for instance, in \cite{lang2002algebra}), we need only show that the functors $\mathcal{H}^i$ are erasable for $i>0$. This is done by induction and a use of the function ring $\mathcal{O}_K$. We claim that $\mathcal{H}^i(\mathfrak{k},\mathcal{O}_K)$ is zero for $i>0$. When $\mathfrak{k}=\mathbb{C}$ and $K=\mathbb{C}$, $\mathcal{O}_K=\mathbb{C}[x]$ and the action of $\mathfrak{k}$ is given by taking derivatives. Thus $\mathcal{H}^1(\mathbb{C},\mathbb{C}[x])=0$ since taking derivative is a surjective map.

Now for general $\mathfrak{k}$, by nilpotency, we have a short exact sequence of Lie algebras $0\to \mathfrak{h}\to \mathfrak{k}\to\mathbb{C}\to 0$. This must split since $\mathbb{C}$ is one dimensional and so we have a covering map $H\rtimes \mathbb{C}\to K$ where $H$ is simply connected. By assumption $K$ is simply connected so the map is an isomorphism. Thus we have an exact sequence of Lie groups $0\to H\to K\to \mathbb{C}\to 0$. Let us consider $\mathcal{H}^*(\mathfrak{k},\mathcal{O}_K)$. By Hochschild-Serre spectral sequence \cite{hochschild1953cohomology}, there is a spectral sequence whose second term is given by $E_2^{*,*}=\mathcal{H}^*(\mathbb{C},\mathcal{H}^*(\mathfrak{h},\mathcal{O}_K))$, that converges to $E_{\infty}^*=\mathcal{H}^{*}(\mathfrak{k},\mathcal{O}_K)$. Since $\mathbb{C}$ is one dimensional, $E_2$ is supported on two columns, the spectral sequence terminates and $\mathcal{H}^{n}(\mathfrak{k},\mathcal{O}_K)=\oplus_{p+q=n}\mathcal{H}^p(\mathbb{C},\mathcal{H}^q(\mathfrak{h},\mathcal{O}_K))$. Consider $\mathcal{H}^q(\mathfrak{h},\mathcal{O}_K)$, we need to understand the module structure of $\mathcal{O}_K$ as an $H$ module. From the isomorphism $H\rtimes \mathbb{C}\cong K$ of Lie groups, we see that there is an isomorphism of algebras \begin{equation}\label{eqLemGnilp}\mathcal{O}_K=\mathcal{O}_H\otimes \mathbb{C}[x],\end{equation} which is described by the following: for $g\in K$, we write $g=h_gc_g$ with $h_g\in H$ and $c_g\in \mathbb{C}$, then the map is given by mapping function $f$ on $K$ to $f(h_gc_g)$ on $H\rtimes \mathbb{C}$. Now to understand the module structure,  if we take an object $f_1\otimes f_2$ where $f_1\in \mathcal{O}_H$ and $f_2\in \mathbb{C}[x]$, for any $h\in H$, $h(f_1\otimes f_2)(g)=f_1\otimes f_2(h^{-1}g)=f_1\otimes f_2(h^{-1}h_gc_g)=f_1(h^{-1}h_g)\otimes f_2(c_g)=((hf_1)\otimes f_2) (g)$. All the equations use the fact that the decomposition of $g=h_gc_g$ is unique. Thus under the above isomorphism \eqref{eqLemGnilp} , $\mathcal{O}_K$ as an $H$ module is nothing but a direct sum of $\mathcal{O}_H$, hence $\mathcal{H}^q(\mathfrak{h},\mathcal{O}_K)=0$ for $q>0$, and $\mathcal{H}^0(\mathfrak{h},\mathcal{O}_K)=\mathcal{O}_K^\mathfrak{h}$, the invariant part of $\mathcal{O}_K$. Again from the identification  \eqref{eqLemGnilp} this is  isomorphic to $\mathbb{C}[x]$. But what is the module structure? Let $f=f_1\otimes f_2$ where $f_1$ is $H$ invariant(it is a constant function in this case), let $c\in \mathbb{C}$, then $cf(g)=f(c^{-1}h_gc_g)=f(c^{-1}hc c^{-1}c_g)$, now since $H$ is a normal subgroup($\mathfrak{h}$ is an ideal), $c^{-1}hc\in H$, and so by the uniqueness of the above decomposition, $cf(g)=f_1(c^{-1}hc)f_2(c^{-1}c_g)=f_1(h)(cf_2)(c_g)$, where we used that $f_1$ is a constant function on $H$. Thus the action on $\mathbb{C}[x]$ is taking derivative and we already see that the cohomology is zero for positive degree. This completes the inductive hypothesis.

Since every $K$ module has an injective resolution by $\CO_K$, we conclude that $\mathcal{H}^i$ are indeed erasable for $i>0$. 

\end{proof}

\subsubsection{Formal Completion}

We are left with computing $\mathrm{End}^{\mathrm{Rep}(\tilde{G}_\CO)}(\CO_{[e]/ \tilde{G}_\CO})$ as an algebraic representation of $\tilde{G}_\CO$. Before going into any details, we would like to comment that this computation may seem complicated, but it is rooted on this simple observation: if $R$ is smooth and $I$ is a complete intersection ideal, then $R/I$ is quasi-isomorphic to its Koszul resolution, and $\mathrm{End}_{R-\mathrm{Mod}}(R/I)$ is an exterior algebra over $R/I$ generated by $(I/I^2)^*$. This is not quite obvious when we replace $R$ by a formally smooth indscheme, since each of the strata may be very singular. In this section, we will introduce formal completion introduced in \cite{gaitsgory2014dg}  to render the situation amenable. 

Let $\mathcal{X}$ be a prestack; then its de-Rham stack is defined by:
\begin{equation}
\mathcal{X}_{dR}(S)=\mathcal{X}(S_{red})
\end{equation} 
where $S_{red}$ is the reduced scheme of $S$. Given a morphism of prestacks $\mathcal{X}\to \mathcal{Y}$, the formal completion is defined by(\cite[Section 6.1]{gaitsgory2014dg}):
\begin{equation}
\widehat{\mathcal{Y}_{\mathcal X}}:=\mathcal{Y}\times_{\mathcal{Y}_{dR}}\mathcal{X}_{dR}.
\end{equation}
This operation behaves well with filtered colimit as explained in \cite[6.1.3]{gaitsgory2014dg}: if $\mathcal{X}=\varinjlim \mathcal{X}_n$ and $\mathcal{Y}=\varinjlim \mathcal{Y}_n$ such that the map $\mathcal{X}\to \mathcal{Y}$ comes from a system of maps $\mathcal{X}_n\to \mathcal{Y}_n$, then: 
\begin{equation}
\widehat{\mathcal{Y}_{\mathcal{X}}}=\varinjlim\widehat{\mathcal{Y}_{n \mathcal{X}_n}}.
\end{equation}

Now assume that $X$ is a locally almost finite type DG scheme and $Y$ an almost finite type DG indscheme, and an embedding $i: X\to Y$, then by \cite[Proposition 6.3.1]{gaitsgory2014dg}, $\widehat{Y_X}$ is a DG indscheme. More-over, from the above we see that:
\begin{equation}
\widehat{Y_X}\cong \varinjlim \widehat{Y_{n X}},
\end{equation}
which in particular means that:
\begin{equation}
\IndCoh(\widehat{Y_X})\cong \varinjlim \IndCoh (\widehat{Y_{n X}}).
\end{equation}
Denote by $\widehat{i}$ the embedding $\widehat{Y_X}\to Y$, and by $\widehat{i_n}$ the embedding of $\widehat{Y_{n X}}\to Y_n$, then by \cite[Proposition 7.4.5]{gaitsgory2014dg}, the adjunction $\mathrm{Id}\to \widehat{i_n}^!\widehat{i_n}_{*,\IndCoh}$ is an equivalence. Taking colimit, we see that $\mathrm{Id}\to \widehat{i}^!\widehat{i}_{*,\IndCoh}$ is an equivalence. If we now consider the sequence of maps:
\begin{equation}
\begin{tikzcd}
X \rar{j} & \widehat{Y_X} \rar{\widehat{i}} & Y,
\end{tikzcd}
\end{equation}
then $i=\widehat{i}\circ j$, and so we have an equivalence of continuous endo-functors of $\IndCoh(X)$:
\begin{equation}
i^!i_{*,\IndCoh}\cong j^!j_{*,\IndCoh}.
\end{equation}

Now let us take $Y=\Gr_G$ and $X$ a miniscule orbit, denote by $\mathcal{X}=X/\tilde{G}_\CO$ and $\mathcal{Y}=Y/\tilde{G}_\CO$. The formal completion $\widehat{Y_X}$ is a DG-indscheme with an $\tilde{G}_\CO$ action, we denote by $\widehat{\mathcal{Y}}_{\mathcal X}$ the quotient stack $\widehat{Y_X}/\tilde{G}_\CO$. We have the following diagram of maps:
\begin{equation}
\begin{tikzcd}
X \rar{j} \arrow[d, "p"] & \widehat{Y_X}\rar{\widehat{i}} \arrow[d, "p"] & Y\arrow[d, "p"]\\
\mathcal X \rar{\overline{j}} & \widehat{\mathcal{Y}}_{\mathcal X}\rar{\overline{\widehat{i}}} & \mathcal{Y} 
\end{tikzcd}
\end{equation}
We comment that the $*$-pushforward and $!$-pullback functors for the category $\IndCoh$ are well defined for the maps $\overline{j}$ and $\overline{\widehat{i}}$. Indeed, one may represent the category $\IndCoh$  over the quotient stacks using the cosimplicial presentation as in equation \eqref{eqcosimplicial}, then the $*$-pushforward and $!$-pullback functors can be defined as the corresponding functors between the cosimplicial categories in equation \eqref{eqcosimplicialcat}. \cite{raskin2020homological} Lemma 6.16.1 and Lemma 6.17.2 together guarantee that this system of functors behave well with the pull-back functors of the cosimplicial categories. 

\newtheorem{pushpull}[LemHequiv]{Lemma}

\begin{pushpull}\label{pushpulllemma}
There is an equivalence of continuous endo-functors of $\IndCoh(\mathcal X)$:
\begin{equation}
\overline{i}^!\overline{i}_{*,\IndCoh}\cong \overline{j}^!\overline{j}_{*,\IndCoh}.
\end{equation}

\end{pushpull}

\begin{proof}
Since $p$ is conservative and t-exact, we need only show that:
\begin{equation}
p^*\overline{i}^!\overline{i}_{*,\IndCoh}\cong p^*\overline{j}^!\overline{j}_{*,\IndCoh}.
\end{equation}
By definition of $\IndCoh^*$ as well as the definition of functors involved, we have $p^*\overline{i}^!\overline{i}_{*,\IndCoh}\cong i^! i_{*,\IndCoh}p^*$, as well as $p^*\overline{j}^!\overline{j}_{*,\IndCoh}=j^!j_{*,\IndCoh}p^*$. These two functors are equivalent as seen from the above discussion. This completes the proof. 
\end{proof}

Recall that we would like to compute $\mathrm{End}^{\mathrm{Rep}(H)}(\overline{i}_{*, \IndCoh}\CO_{\mathcal{X}})$. By adjunction:
\begin{equation}
\mathrm{End}^{\mathrm{Rep}(H)}(\overline{i}_{*, \IndCoh}\CO_{\mathcal{X}})\cong \mathrm{Hom}^{\mathrm{Rep}(H)}(\CO_{\mathcal{X}}, \overline{i}^!\overline{i}_{*, \IndCoh}\CO_{\mathcal{X}}).
\end{equation}
By Lemma \ref{pushpulllemma} we have:
\begin{equation}
\mathrm{Hom}^{\mathrm{Rep}(H)}(\CO_{\mathcal{X}}, \overline{i}^!\overline{i}_{*, \IndCoh}\CO_{\mathcal{X}})\cong \mathrm{Hom}^{\mathrm{Rep}(H)}(\CO_{\mathcal{X}}, \overline{j}^!\overline{j}_{*, \IndCoh}\CO_{\mathcal{X}}).
\end{equation}
 Thus we have transfered the computation onto the formal completion. In the next section, we will specialize to the case when $X=[e]$ and $Y=\Gr_G$, and explicitly understand this formal completion using the idea of formal geometry studied in  \cite{gaitsgory2020study}.
 
 \subsubsection{Formal Groups and Lie Algebras}
 
 In  \cite[Chapter 7]{gaitsgory2020study}, the authors studied formal groups, and showed that the category of formal groups over a prestack $\mathcal X$ is equivalent to that of Lie algebra objects in $\IndCoh(\mathcal X)$. Let us briefly recall the important notions here. Let $\mathcal X$ be a locally almost finite type prestack (see  \cite{gaitsgory2020study} for what it means). Denote by $\mathrm{FormMod}_{/ \mathcal X}$ the category of locally almost finite type stacks $\mathcal Z$ over $\mathcal X$ such that the map $\mathcal{Z}\to \mathcal{X}$ is inf-schematic and induces an equivalence $\mathcal{Z}_{red}\cong \mathcal{X}_{red}$ (\cite[Chapter 5, 1.1.1]{gaitsgory2020study}).  A formal group over $\mathcal X$ is a group object in $\mathrm{FormMod}_{/ \mathcal X}$. This category is denoted by $\mathrm{Grp}(\mathrm{FormMod}_{/ \mathcal X})$. On the other hand, consider the category of Lie algebra objects in $\IndCoh(\mathcal X)$, which we denote by $\mathrm{LieAlg}(\IndCoh(\mathcal X))$. The result of \cite[Chapter 7]{gaitsgory2020study}, more specifically Theorem 3.1.4, states that there is an equivalence:
 \begin{equation}
 \mathrm{Grp}(\mathrm{FormMod}_{/ \mathcal X})\cong \mathrm{LieAlg}(\IndCoh(\mathcal X)). 
 \end{equation}
 The idea of this is that given a formal group $\mathcal Y$ over $\mathcal X$, the object $\pi_{*,\IndCoh}(\omega_{\mathcal Y})$, the pushforward of the dualizing sheaf, has the structure of a cocommutative Hopf algebra. This is the universal enveloping algebra of the Lie algebra associated to $\mathcal Y$. 
 
When $\mathcal{X}=\mathrm{pt}$, then the category $\mathrm{LieAlg}(\IndCoh(\mathrm{pt}))$ is the category of DG Lie algebras in $\mathrm{Vect}$ studied in \cite{lurie2011derived}. In the special case when $\mathfrak{g}$ is a Lie algebra concentrated in degree $0$, the formal moduli problem is simply $\widehat{\mathfrak{g}_0}$, the formal completion of $\mathfrak{g}$ at $0$(\cite[Construction 2.2.13.]{lurie2011derived}). The formal group structure is given by the Baker–Campbell–Hausdorff formula. 

Let us now apply this to the case when $X=[e]$ and $Y=\Gr_G$, we have:

\newtheorem{formalgroup}[LemHequiv]{Lemma}

\begin{formalgroup}\label{formalgroup}
The formal completion $\widehat{Y_X}$ is a formal group whose Lie algebra is $z^{-1}\mathfrak{g}[z^{-1}]$.

\end{formalgroup} 
 
 \begin{proof}
 From the discussion of Section \ref{BFNspace}, the group ind-scheme $L^{<0}G$ is an open neighborhood of $X$ in $Y$, and so $\widehat{Y_X}\cong \widehat{L^{<0}G_{X}}$. Now $\widehat{L^{<0}G_{X}}$ is a formal group whose associated Lie algebra is $z^{-1}\mathfrak{g}[z^{-1}]$.

 \end{proof}
 
 Denote by $L^{<0}\mathfrak{g}$ the Lie algebra of $L^{<0}G$. By \cite[Chapter 7, Theorem 3.1.4]{gaitsgory2020study}, we see that  $\widehat{Y_X}$ is equivalent to $\widehat{L^{<0}\mathfrak{g}_{0}}$, the formal completion of $L^{<0}\mathfrak{g}$ at $0$. The action of $\tilde{G}_\CO$ is given by conjugation on $L^{<0}\mathfrak{g}\cong \mathfrak{g}(\CK)/\mathfrak{g}(\CO)$. Again denote by $\mathcal{X}=X/\tilde{G}_\CO$ and $\mathcal{Y}=Y/\tilde{G}_\CO$. Recall the morphism $\overline{j}: \mathcal{X}\to \widehat{\mathcal Y}_{\mathcal X}$ and $\overline{i}: \mathcal{X}\to\mathcal{Y}$. We claim:
 
 \newtheorem{pushpullid}[LemHequiv]{Proposition}
 
 \begin{pushpullid}\label{pushpullid}
 There is an equivalence of continuous endofunctors on $\IndCoh(\mathcal{X})$
 \begin{equation}\label{eqpushpullid}
 \overline{i}^!\overline{i}_{*,\IndCoh}\cong \mathrm{Sym}^{\bullet}\!\!\left(L^{<0}\mathfrak{g}[-1]\right)\otimes -
 \end{equation}
 where $L^{<0}\mathfrak{g}$ is understood as a $G(\CO)$ module under conjugation action.
 \end{pushpullid}
 
 \begin{proof}
 By Lemma \ref{pushpulllemma}, we can replace the left hand side of equation \eqref{eqpushpullid} by $ \overline{j}^!\overline{j}_{*,\IndCoh}$. Consider the following diagram:
 \begin{equation}
 \begin{tikzcd}
 \mathcal{X}\rar{\overline{j}} & \tilde{G}_{\CO}\setminus \widehat{L^{<0}\mathfrak{g}_0}\rar & G(\CO)\setminus  L^{<0}\mathfrak{g}
 \end{tikzcd}
 \end{equation}
 Denote by $\overline{i}_g$ the inclusion $\mathcal{X}\to G(\CO)\setminus  L^{<0}\mathfrak{g}$, Lemma \ref{pushpulllemma} again implies:
 \begin{equation}
  \overline{j}^!\overline{j}_{*,\IndCoh}\cong \overline{i}_g^!\overline{i}_{g,*,\IndCoh}.
 \end{equation}
 The latter can be computed explicitly using a Koszul resolution, and the result follows. 
 
 \end{proof}
 
 We can now prove:

\newtheorem{CorExtidentity}[LemHequiv]{Corollary}

\begin{CorExtidentity}\label{CorExtidentity}
There is a quasi-isomorphism of $\tilde{G}_\CO$ vector spaces:
\begin{equation}\label{eqExtidentity}
\mathrm{End}^{\mathrm{Rep}(\tilde{G}_\CO)}(\CO_{[e]/\tilde{G}_\CO})\cong \mathrm{Sym}^{\bullet}\!\!\left(\mathfrak{g}(\CK)/\mathfrak{g}(\CO)[-1]\right).
\end{equation}
\end{CorExtidentity}

\begin{proof}
By Proposition \ref{pushpullid}:
\begin{equation}\label{corextidentity}
\mathrm{End}^{\mathrm{Rep}(\tilde{G}_\CO)}(\CO_{[e]/\tilde{G}_\CO})\cong\mathrm{Hom}^{\mathrm{Rep}(\tilde{G}_\CO)}\left(\CO_{[e]/\tilde{G}_\CO}, \mathrm{Sym}^{\bullet}\!\!\left(\mathfrak{g}(\CK)/\mathfrak{g}(\CO)[-1]\right) \otimes \CO_{[e]/\tilde{G}_\CO} \right).
\end{equation}
Since $\CO_{[e]/\tilde{G}_\CO}$ is simply the trivial representation of $\tilde{G}_\CO$, the right hand side of equation \eqref{corextidentity} can be identified as the right hand side of equation \eqref{eqExtidentity}. This completes the proof. 
\end{proof}

 \subsubsection{Other Miniscule Orbits}

We can in fact use this technique for other miniscule orbits of $\Gr_G$. Let us now take $X$ to be a miniscule orbit and $Y=\Gr_G$. Denote by $\mathcal X$ and $\mathcal Y$ the quotients of $X$ and $Y$ by $\tilde{G}_\CO$. Choose $[g]$ a point in $\mathcal X$, let $\tilde{P}$ be the stabilizer of $[g]$ in $\tilde{G}_\CO$, then the there is an equivalence of prestacks:
\begin{equation}
\mathcal{X}\cong \tilde{G}_\CO\setminus \tilde{G}_\CO/\tilde{P}\cong \mathbb{B}\tilde{P}.
\end{equation}
Under this, the map $\overline{i}: \mathcal{X}\to \mathcal{Y}$ corresponds to the map of schemes:
\begin{equation}\label{eqGOP}
\begin{tikzcd}
\tilde{P}\setminus\mathrm{pt}\rar{j} & \tilde{P}\setminus \mathrm{Gr}_G \rar{m} & \tilde{G}_{\CO}\setminus \mathrm{Gr}_G
\end{tikzcd}
\end{equation}
where the map $j$ is the embedding of $\mathrm{pt}$ as $[g^{-1}]$. Let $V$ be the $\tilde{P}$ module given by:
\begin{equation}
\mathfrak{g}(\CK)/(\mathfrak{g}(\CO)+ g\mathfrak{g}(\CO)g^{-1}).
\end{equation}
We prove:

\newtheorem{minisculeorbit}[LemHequiv]{Proposition}

\begin{minisculeorbit}\label{minisculeorbit}
There is a quasi-isomorphism of objects in $\IndCoh(\mathcal{X})$:
\begin{equation}
\overline{i}^!\overline{i}_{*,\IndCoh}(\CO_{\mathcal X})\cong \tilde{G}_\CO\times_{\tilde{P}} \mathrm{Sym}^{\bullet}\!(V[-1]).
\end{equation}

\end{minisculeorbit}

\begin{proof}
Using the presentation of $\mathcal X$ in equation \eqref{eqGOP}, we would like to show that:
 \begin{equation}
 j^!m^!m_{*,\IndCoh}j_{*,\IndCoh}(\CO_{\mathcal X})= \mathrm{Sym}^{\bullet}\!(V[-1])
 \end{equation}
as a module of $\tilde{P}$. Let us understand the composition $m^!m_{*,\IndCoh}j_{*,\IndCoh}(\CO_{\mathcal X})$, consider the following Cartesian diagram:
\begin{equation}
\begin{tikzcd}
\tilde{P}\setminus \tilde{G}_\CO / \tilde{P} \rar{\tilde{m}} \arrow[d, "\tilde{\overline{i}}"] & \tilde{P}\setminus\mathrm{pt} \arrow[d, "\overline{i}"]\\
\tilde{P}\setminus \mathrm{Gr}_G \rar{m} & \tilde{G}_{\CO}\setminus \mathrm{Gr}_G
\end{tikzcd}
\end{equation}
Here $\tilde{m}$ is the projection of $\tilde{G}_\CO/\tilde{P}$ to a point, and $\tilde{\overline{i}}$ is induced by the embedding $X\to \Gr_G$. 
By base-change property \cite[Proposition 2.9.2]{gaitsgory2014dg}, we have: 
\begin{equation}
m^!\overline{i}_{*,\IndCoh}\cong \tilde{\overline{i}}_{*,\IndCoh}\tilde{m}^!.
\end{equation}
Thus the object $m^!\overline{i}_{*,\IndCoh}(\CO_{\mathcal X})$ is $\tilde{\overline{i}}_{*,\IndCoh}\omega_{X}$, where $\omega_X$ is the dualizing sheaf of $X$ with the canonical $\tilde{P}-$equivariant structure. In our case, since $X=\tilde{G}_\CO/\tilde{P}$, $\omega_X$ is the line bundle over $X$ associated to the  one dimensional $\tilde{P}$ representation:
\begin{equation}
L_{\mathrm{top}}=\mathrm{Sym}^{\mathrm{top}}\!\left(\mathfrak{g}_\CO/\mathfrak{p}[1]\right)
\end{equation}
Here $\mathfrak{g}_\CO/\mathfrak{p}[1]$ is a finite dimensional vector space in cohomological degree $-1$, and so the exterior algebra has finite cohomological degree. The representation $L_{\mathrm{top}}$ is the top degree part of the exterior algebra, and is in cohomological degree $-\mathrm{dim}(X)$. Let us now employ the idea of formal completion. Consider the Cartesian diagram:
\begin{equation}
\begin{tikzcd}
\tilde{P}\setminus \widehat{X_{[g]}} \rar \arrow[d] & \tilde{P}\setminus \tilde{G}_\CO / \tilde{P}\arrow[d]\\
\tilde{P}\setminus \widehat{Y_{[g]}} \rar & \tilde{P}\setminus \mathrm{Gr}_G
\end{tikzcd}
\end{equation}
By base-change property \cite[Proposition 2.9.2]{gaitsgory2014dg}, the shriek pullback of $\tilde{\overline{i}}_{*,\IndCoh}\omega_{X}$ to $\widehat{Y_{[g]}}$ is the pushforward of the dualizing sheaf of $\widehat{X_{[g]}}$ to $\widehat{Y_{[g]}}$. The advantage is that these local completions have very explicit descriptions. Indeed, by \cite[Chapter 7, Theorem 3.1.4]{gaitsgory2020study}, the space $\widehat{X_{[g]}}$ is equivalent to the completion of $\mathfrak{g}_{\CO}/\mathfrak{p}$ at $0$. Similarly, the space $\widehat{Y_{[g]}}$ is equivalent to the completion of $L^{<0}\mathfrak{g}$ at $0$. The map $\widehat{X_{[g]}}\to \widehat{Y_{[g]}}$ corresponds to the embedding of the following $\tilde{P}$ modules:
\begin{equation}\label{embedlie}
\psi: \mathfrak{g}_{\CO}/\mathfrak{p}\to \mathfrak{g}(\CK)/\mathfrak{g}(\CO), ~ H\to gHg^{-1}. 
\end{equation}
We can thus transfer to the following diagram:
\begin{equation}
\begin{tikzcd}
 & \tilde{P}\setminus \widehat{X_{[g]}}\rar{\widehat{\varphi}} \arrow[d,"\widehat{\psi}"] & \tilde{P}\setminus \mathfrak{g}_{\CO}/\mathfrak{p} \arrow[d,"\psi"]\\
 \tilde{P}\setminus \mathrm{pt} \rar{\widehat{j}} &  \tilde{P}\setminus\widehat{Y_{[g]}} \rar{\widehat{\phi}} & \tilde{P}\setminus \mathfrak{g}(\CK)/\mathfrak{g}(\CO)
\end{tikzcd}
\end{equation}
with which we can derive: 
\begin{equation}
j^!\tilde{\overline{i}}_{*,\IndCoh}\omega_{X}\cong \widehat{j}^! \widehat{\phi}^! \psi_{*,\IndCoh} (\omega_{ \mathfrak{g}_{\CO}/\mathfrak{p}}).
\end{equation}
Here the sheaf $\omega_{ \mathfrak{g}_{\CO}/\mathfrak{p}}$ is the structure sheaf of $\mathfrak{g}_{\CO}/\mathfrak{p}$ tensored with the representation $L_{\mathrm{top}}$. We now have:
\begin{equation}
\widehat{j}^! \widehat{\phi}^! \psi_{*,\IndCoh}(\omega_{ \mathfrak{g}_{\CO}/\mathfrak{p}})\cong \mathrm{Hom}^{\mathbb{B}\tilde{P}}\left((\widehat{j}\circ\widehat{\phi})_{*,\IndCoh}(\CO_{\mathbb{B}\tilde{P}}), \psi_{*,\IndCoh} (\omega_{ \mathfrak{g}_{\CO}/\mathfrak{p}})\right).
\end{equation}
The right hand side can be computed using a Koszul resolution of $(\widehat{j}\circ\widehat{\phi})_{*,\IndCoh}(\CO_{\mathbb{B}\tilde{P}})$, and the result is the following complex:
\begin{equation}\label{eqKoszulmini}
\mathrm{Sym}^{\bullet}\!\!\left(\mathfrak{g}(\CK)/\mathfrak{g}(\CO)[-1]\right)\otimes \mathbb{C}[\mathfrak{g}_\CO/\mathfrak{p}]\otimes L_{\mathrm{top}},
\end{equation}
together with a differential induced from the Koszul resolution. Here $\mathbb{C}[\mathfrak{g}_\CO/\mathfrak{p}]$ denotes the algebra of functions on $\mathfrak{g}_\CO/\mathfrak{p}$. The nonzero part of the Koszul differential lies in: 
\begin{equation}
\mathrm{Sym}^{\bullet}\!\!\left(\mathfrak{g}(\CO)/\mathfrak{p}[-1]\right) \otimes \mathbb{C}[\mathfrak{g}_\CO/\mathfrak{p}] \otimes L_{\mathrm{top}}\cong \mathrm{Sym}^{\bullet}\!\!\left(\left( \mathfrak{g}(\CO)/\mathfrak{p}\right)^*[1]\right)\otimes  \mathbb{C}[\mathfrak{g}_\CO/\mathfrak{p}] .
\end{equation}
The quasi-isomorphism is due to tensoring with $L_{\mathrm{top}}$, which makes this into a usual Koszul complex. The cohomology of this complex is $\mathbb{C}$ in degree $0$, and so the cohomology of the complex in equation \eqref{eqKoszulmini} is thus identified with $\mathrm{Sym}^{\bullet}\!(V[-1])$. This completes the proof. 
\end{proof}

\subsubsection{Conclusion}

Using Proposition \ref{PropGOequiv} and Corollary \ref{CorExtidentity}, we obtain:

\newtheorem{ExtidentityThm}[LemHequiv]{Theorem}

\begin{ExtidentityThm}\label{Extidentitypure}
There is a quasi-isomorphism of $\mathbb{C}^*$ vector spaces:
\begin{equation}
\mathrm{End}^{\mathbb{B}\mathbb{C}^*}(\CO_{[e]/\tilde{G}_\CO})\cong \left[ \mathrm{Sym}^{\bullet}\!\!\left(\left(\mathfrak{g}(\CK)/\mathfrak{g}(\CO)\oplus \left( z\mathfrak{g}(\CO)\right)^*\right)[-1]\right)\right]^G,
\end{equation}
where $[-]^G$ is taking ordinary $G$ invariants. This space coincides with $\pi_0\mathcal{V}_{G,0}$ of equation \eqref{eqpoissondualpure} after shifting the degree of $\mathfrak{g}(\CK)/\mathfrak{g}(\CO)$ to $-1$.  

\end{ExtidentityThm}

\noindent\textbf{Remark}. As remarked in \cite{oh2020poisson}, the character of the above space is given by:
\begin{equation}
\frac{1}{|W|}\oint_T \frac{\mathrm{d}s}{2\pi i s}\prod\limits_{\alpha\text{ roots}}  (1-s^\alpha) \left[(q)_{\infty}^{2\mathrm{rank}(G)}\prod\limits_{\alpha\text{ roots}}(qs^\alpha; q)_{\infty}^2\right],
\end{equation}
which reproduces Schur index of a pure gauge theory.

\subsection{The Abelian Gauge Group Case}\label{AbelianGaugeGroup}

As a non-trivial but illuminating example, we consider the case $T=\mathbb{C}^*$. The space $\Gr_{T}$ is the moduli space of lattices in $\CK$. This is an ind-projective variety, but is not reduced. There is a closed embedding: $\mathrm{Gr}_{T,red}\to \mathrm{Gr}_T$, where $\Gr_{T,red}$ is the reduced scheme of $\Gr_T$, and is a discrete set of points labeled by $\mathbb{Z}$. In particular, the identity coset $e$ is a closed reduced point in $\mathrm{Gr}_T$. The Schur index for the abelian gauge group $T$ can be reproduced by the endomorphism algebra:
\begin{equation}
\mathrm{End}^{\mathbb{B}\mathbb{C}^*}(\CO_{[e]/T(\CO) \rtimes \mathbb{C}^*}).
\end{equation}
The argument from Section \ref{PuregaugessG} can still be applied. The relative tangent complex of $e\to \Gr_T$  is the space $\CK/\CO[1]$, and the action of $\CO^*$ is trivial. We have the following quasi-isomorphism:
\begin{equation}
\mathrm{End}^{\mathbb{B}\mathbb{C}^*}(\CO_{\mathbb{B}(T(\CO) \rtimes \mathbb{C}^*)})\cong \mathrm{Sym}^{\bullet}\!\left(\left(\CK/\CO\oplus (z\CO)^*\right)[-1]\right).
\end{equation}
In this case, the Poisson algebra $\mathcal{V}_{G,0}$ is just $bc$ ghost with trivial differential, and the space of local operators $\pi_0\mathcal{V}_{G,0}$ coincides with the above as a vector space.

It is important in this example that we use the non-reduced scheme. The embedding $i:\mathrm{Gr}_{T,red}\to \mathrm{Gr}_T$  induces an isomorphism on $K_0$ group,  but the functors $i_{*,\IndCoh}$ and $i^{!}$ are not inverses to each other. In fact, Proposition \ref{pushpullid} implies that $i^{!}i_{*,\IndCoh}(\CO_{z^n})=\mathrm{Sym}^{\bullet}\!\left(\CK/z^n\CO[-1]\right)$, and this factor would be absent if one used $\Gr_{T,red}$ istead of $\Gr_T$.  

In fact, when $T=\mathbb{C}^*$, the affine Grassmannian $\Gr_T$ has a rather simple form: let $\widehat{W}$ be the functor assigning to an algebra $R$ the space of invertible elements in $1+z^{-1}R[z^{-1}]$; then one has a decomposition:
\begin{equation}
\Gr_T=\widehat{W}\times \Gr_{T,red}.
\end{equation}
Here $\widehat{W}$ is a formal group whose reduced scheme is a point. It is called the formal group of Witt vectors. Its Lie algebra is $L^{<0}\mathbb{C}$. We thus have an equivalence of formal stacks:
\begin{equation}
\widehat{W}=\widehat{L^{<0}\mathbb{C}_0},
\end{equation}
from which we conclude that $\Gr_T=\widehat{L^{<0}\mathbb{C}_0}\times \mathbb{Z}$. In conclusion, when $T=\mathbb{C}^*$, even without taking formal completion, the affine Grassmannian $\Gr_T$ decomposes into a product of $\mathbb{Z}$ with the formal completion of a vector space at a point. Consequently, one can compute the space of local operators at the junction of t'Hooft  line defects labeled by $\Gr_{T,red}$, by reducing the computation to $\widehat{W}$. In principle, as far as the index is concerned, one may consider more general sheaves over $\widehat{W}$, by writing out the composition series of a coherent sheaf using the structure sheaf of the reduced point. We will not go in this direction here.

\subsection{Gauge Theory with Matter}\label{gaugetheorywithmatter}

In this section, we will generalize the computation of Section \ref{PuregaugessG} to gauge theories with matter. Let $V$ be a representation of $G$. Recall that the BFN space is defined by the base change diagram:
\begin{equation}
\begin{tikzcd}
\mathcal{R}_{G,V} \rar \arrow[d] &  V(\mathcal{O})\arrow[d] \\  G(\mathcal{K})\times_{G(\mathcal{O})}V(\mathcal{O}) \rar &  V(\mathcal{K})
\end{tikzcd}
\end{equation} 
We add to this another base-change diagram:
\begin{equation}
\begin{tikzcd}
Z \rar \arrow[d] & \mathcal{R}_{G,V} \rar \arrow[d] &  V(\mathcal{O})\arrow[d] \\ e\times V(\CO)\rar & G(\mathcal{K})\times_{G(\mathcal{O})}V(\mathcal{O}) \rar &  V(\mathcal{K})\\
\end{tikzcd}
\end{equation} 
Here $Z=V(\CO)\times_{V(\CK)}V(\CO)$ can be described as $V(\CO)\times V(\CK)/V(\CO)[-1]$. The identity line is the pushforward of structure sheaf of $V(\CO)$ along the embedding $i:V(\CO)\to \mathcal{R}_{G,N}$. Note that this is a classical scheme embedded into a derived scheme. We will label the maps:
\begin{equation}
\begin{tikzcd}
V(\CO) \rar{l} &  Z \rar{m} \arrow[d,"p_1"] & \mathcal{R}_{G,V} \rar \arrow[d, "p_2"] &  V(\mathcal{O})\arrow[d] \\ &  e\times V(\CO)\rar{j} & G(\mathcal{K})\times_{G(\mathcal{O})}V(\mathcal{O}) \rar &  V(\mathcal{K})
\end{tikzcd}
\end{equation} 
The Schur index is then the graded Euler character of:
\begin{equation}
\mathrm{End}^{\mathbb{B}\mathbb{C}^*}\left( i_{*,\IndCoh}(\CO_{V(\CO)/\tilde{G}_\CO})\right).
\end{equation}

In the following, we will write $\overline{X}$ for the quotient stack $X/\tilde{G}_\CO$, in order to avoid clustering of notations. We will also omit the $\IndCoh$ for all the push-forward functors. To do this computation, fix again an ind-scheme structure of $\Gr_G$ and $\CBFN$ similar to Section \ref{BFNspace}. We make use of the following diagram:
\begin{equation}
\begin{tikzcd}
V(\CO) \rar{l_n} &  Z_n \rar{m_n} \arrow[d,"p_1"] & \mathcal{R}_{G,V,n} \rar \arrow[d, "p_2"] &  V(\mathcal{O})\arrow[d] \\ &  e\times V(\CO)\rar{j_n} & G(\mathcal{K})_n\times_{G(\mathcal{O})}V(\mathcal{O}) \rar &  z^{-N}V(\mathcal{O})
\end{tikzcd}
\end{equation} 
such that $m_n\circ l_n=i_n$. Since $\CBFN$ is a colimit of $ \mathcal{R}_{G,V,n}$, by equation \eqref{eqcolimTreal}:
\begin{equation}
\mathrm{End}^{\mathbb{B}\mathbb{C}^*}\left( i_{*}(\CO_{\overline{V(\CO)}})\right)=\varinjlim\limits_{n} \mathrm{End}^{\mathbb{B}\mathbb{C}^*}\left( i_{n,*}(\CO_{\overline{V(\CO)}})\right).
\end{equation}
Lemma \ref{LemHRep} implies that $\mathrm{End}^{\mathbb{B}\mathbb{C}^*}\left(  i_{n,*} (\CO_{\overline{V(\CO)}})\right)$ this is the $G(\CO)$ invariants of:
\begin{equation}\label{EndBFNnogroup}
\mathrm{End}^{\mathrm{Rep}(\tilde{G}_\CO)}\left(  i_{n,*}(\CO_{\overline{V(\CO)}})\right).
\end{equation}
Let us compute this vector space using adjunctions. Since $i_n=m_n\circ l_n$, from the adjunction pair $(m_{n,*},m_n^{!})$, one has:
\begin{equation}
\mathrm{End}^{\mathrm{Rep}(\tilde{G}_\CO)}\left( i_{n,*}(\CO_{\overline{V(\CO)}})\right)\cong \mathrm{Hom}^{\mathrm{Rep}(\tilde{G}_\CO)}\left( l_{n,*}(\CO_{\overline{V(\CO)}}), m_n^!i_{n,*}(\CO_{\overline{V(\CO)}})\right).
\end{equation}
 As $Z_n$ is a very explicit DG scheme with a very explicit action of $\tilde{G}_\CO$, one can write an explicit projective resolution of $l_{n,*}(\CO_{\overline{V(\CO)}})$ given by the Koszul complex:
\begin{equation}\label{eqresolutionNO}
l_{n,*}(\CO_{\overline{V(\CO)}})\cong \CO_{\overline{Z_n}} \otimes \mathrm{Sym}^{\bullet}\!\!\left((z^{-N}V(\CO)/V(\CO))^*[2]\right).
\end{equation}
 together with the differential given by the usual Koszul differential. This is a quasi-isomorphism of $\tilde{G}_\CO$ equivariant sheaves. Substituting the resolution of equation \eqref{eqresolutionNO} into the above equation, one has:
\begin{equation}
\begin{aligned}
\mathrm{Hom}&^{\mathrm{Rep}(\tilde{G}_\CO)}\left(l_{n,*}(\CO_{\overline{V(\CO)}}), m_n^{!}i_{n,*}(\CO_{\overline{V(\CO)}})\right)\\ &\cong \mathrm{Hom}^{\mathrm{Rep}(\tilde{G}_\CO)}\left(\mathcal{O}_{\overline{Z_n}}, m_n^{!}i_{n,*}(\CO_{\overline{V(\CO)}})\right)\otimes \mathrm{Sym}^{\bullet}\!\!\left(z^{-N}V(\CO)/V(\CO)[-2]\right).
\end{aligned}
\end{equation}
By definition, $\mathcal{O}_{\overline{Z_n}}=(p_1)^*\CO_{\overline{V(\CO)}}$, using push-pull adjunction, one has:
\begin{equation}
\mathrm{Hom}^{\mathrm{Rep}(\tilde{G}_\CO)}\left(\mathcal{O}_{\overline{Z_n}}, m_n^{!}i_{n,*}(\CO_{\overline{V(\CO)}})\right)\cong \mathrm{Hom}^{\mathrm{Rep}(\tilde{G}_\CO)}\left(\CO_{\overline{V(\CO)}}, (p_1)_*m_n^{!}i_{n,*}(\CO_{\overline{V(\CO)}})\right).
\end{equation}
 We then apply the base-change property established in \cite{raskin2020homological} Lemma 6.16.1, namely that $(p_1)_*m_n^!\cong j_n^! (p_2)_*$, which implies:
 \begin{equation}
  \mathrm{Hom}^{\mathrm{Rep}(\tilde{G}_\CO)}\left(\CO_{\overline{V(\CO)}}, (p_1)_*m_n^{!}i_{n,*}(\CO_{\overline{V(\CO)}})\right)\cong  \mathrm{Hom}^{\mathrm{Rep}(\tilde{G}_\CO)}\left(\CO_{\overline{V(\CO)}}, j_n^!j_{n,*}(\CO_{\overline{V(\CO)}})\right).
 \end{equation}
 To make contact with the affine Grassmannian, we now consider the following Cartesian diagram:
\begin{equation}
\begin{tikzcd}
e\times V(\CO)\rar{j_n} \arrow[d, "q_1"] & G(\mathcal{K})_{n}\times_{G(\mathcal{O})}V(\CO) \arrow[d,"q_2"]\\e \rar{k_n} &  \Gr_{G,n}.
\end{tikzcd}
\end{equation}
Since $\CO_{\overline{V(\CO)}}\cong q_1^*\CO_{\overline{[e]}}$, by pull-push adjunction:
\begin{equation}
\mathrm{Hom}^{\mathrm{Rep}(\tilde{G}_\CO)}\left(\CO_{\overline{V(\CO)}}, j_n^!j_{n,*}(\CO_{\overline{V(\CO)}})\right)\cong \mathrm{Hom}^{\mathrm{Rep}(\tilde{G}_\CO)}\left(\CO_{\overline{[e]}}, (q_1)_*j_n^!j_{n,*}(\CO_{\overline{V(\CO)}})\right).
\end{equation}
By base-change formula again, $ (q_1)_*j_n^!\cong k_{n}^!(q_2)_*$, we obtain:
\begin{equation}
\mathrm{Hom}^{\mathrm{Rep}(\tilde{G}_\CO)}\left(\CO_{\overline{[e]}}, (q_1)_*j_n^!j_{n,*}(\CO_{\overline{V(\CO)}})\right)\cong \mathrm{End}^{\mathrm{Rep}(\tilde{G}_\CO)}(k_{n,*} \CO_{\overline{[e]}})\otimes \mathbb{C}[V(\CO)].
\end{equation}
 Here $V(\CO)$ is in cohomological degree $0$. By taking the colimit and applying Proposition \ref{colimRepH}, we find that the underlying $\tilde{G}_\CO$ representation of $\mathrm{End}^{\mathrm{Rep}(\tilde{G}_\CO)}( i_{*}(\CO_V))$ can be identified with:
 \begin{equation}
 \mathbb{C}[V(\CO)]\otimes  \mathrm{Sym}\left(V(\CK)/V(\CO)\right)\otimes  \mathrm{End}^{\mathrm{Rep}(\tilde{G}_\CO)}(\CO_{\overline{[e]}}).
 \end{equation}
 Here $\mathrm{End}^{\mathrm{Rep}(\tilde{G}_\CO)}(\CO_{\overline{[e]}})$ is computed in equation \eqref{eqExtidentity}. 

As far as the character is concerned, the above computation thus gives us the desired $\tilde{G}_\CO$ module. However, the differential is kept obscured in this computation. To analyze the differential, we will use formal completion. Recall the following diagram:
\begin{equation}
\begin{tikzcd}
\mathcal{R}_{G,V} \rar \arrow[d] &  V(\mathcal{O})\arrow[d] \\  G(\mathcal{K})\times_{G(\mathcal{O})}V(\mathcal{O}) \rar &  V(\mathcal{K})
\end{tikzcd}
\end{equation}
Denote by $\widehat{\mathcal{T}}$ the formal completion of $G(\mathcal{K})\times_{G(\mathcal{O})}V(\mathcal{O})$ along $[e]\times V(\CO)$. This is a $\tilde{G}_\CO$-equivariant formal scheme over $V(\CO)$. It is clear that it is isomorphic to $\widehat{\Gr_{G,[e]}}\times V(\CO)$ where $\widehat{\Gr_{G,[e]}}$ is the formal completion of $\Gr_G$ along $[e]$. As already discussed in Lemma \ref{formalgroup}, the space $\widehat{\Gr_{G,[e]}}$ is a formal group, and thus by \cite[Theorem 3.1.4]{gaitsgory2020study}, it is isomorphic, as a formal group, to the formal completion of its Lie algebra at $0$, namely $\widehat{L^{<0}\mathfrak{g}_{0}}$. We define $\widehat{\mathcal R}$ by the following diagram:
 \begin{equation}
\begin{tikzcd}
\widehat{\mathcal R}\rar\arrow[d] & \mathcal{R}_{G,V} \rar \arrow[d] &  V(\mathcal{O})\arrow[d] \\ \widehat{L^{<0}\mathfrak{g}_0}\times V(\CO)\rar &  G(\mathcal{K})\times_{G(\mathcal{O})}V(\mathcal{O}) \rar &  V(\mathcal{K})
\end{tikzcd}
\end{equation}
By \cite[Section 6.1.3 (iv)]{gaitsgory2014dg}, $\widehat{\mathcal R}$ can be identified as the formal completion of $\mathcal{R}_{G,V}$ along $V(\CO)\times_{V(\CK)}V(\CO)$. Let $\hat{i}$ be the embedding $V(\CO)\to \widehat{\mathcal R}$, then just as in Lemma \ref{pushpulllemma}, we have:
\begin{equation}
\mathrm{Hom}^{\mathrm{Rep}(\tilde{G}_\CO)}\left(\CO_{\overline{V(\CO)}}, \hat{i}^!\hat{i}_*\CO_{\overline{V(\CO)}}\right)\cong \mathrm{Hom}^{\mathrm{Rep}(\tilde{G}_\CO)}\left(\CO_{\overline{V(\CO)}}, i^!i_*\CO_{\overline{V(\CO)}}\right).
\end{equation}
The advantage of this construction is the following: the space $\widehat{\mathcal R}$ is an explicit DG ind-scheme whose underlying pro-algebra is represented by the pro-algebra of functions on the following ind-scheme:
\begin{equation}
\widehat{L^{<0}\mathfrak{g}_0}\times V(\CO)\times V(\CK)/V(\CO)[-1],
\end{equation} 
and has a differential $D$ as described in Section \ref{BFNspace}, induced by the formal group action. We will denote by $A$ the pro-algebra defining this DG ind-scheme. It is worth writing down this differential explicitly here. Choose a basis $v^i$ for $V$, let $\rho_{i}^j$ be the matrix elements of $\mathfrak{g}$ action on $V$, namely:
\begin{equation}
Xv^j=\sum_i \rho_{i}^j(X)v^i.
\end{equation}
Denote by $\rho_{i,n}^j$ the corresponding linear function on $L^{<0}\mathfrak{g}$, by $v_{i,n}^*$ the corresponding linear functions on $V(\CO)$, and by $w_{i,n}^*$ the linear functions on $V(\CK)/V(\CO)[-1]$. Note that $w_{i,n}^*$ are odd variables. The differential $D$ can be expressed as:
\begin{equation}
Dw_{i,n}^*=\sum_{j, m+k=n}\rho_{i,m}^j\otimes v_{j,k}^*+\frac{1}{2}\sum_{j_1,j_2, m_1+m_2+k=n} \rho_{i,m_1}^{j_1}\rho_{j_1, m_2}^{j_2}\otimes v_{j_2,k}^*+\cdots.
\end{equation}
We comment that this comes from exponentiating the action of $X$, and the first term vanish because functions on $V(\CK)/V(\CO)$ is zero on $V(\CO)$. This differential should be understood in the pro-algebra otherwise the summation would be infinite.

The identity line is the structure sheaf of $e\times V(\CO)$, and we would like to use a Koszul resolution:
\begin{equation}
\mathbb{C}[V(\CO)]\cong \mathrm{Sym}^{\bullet}\!\!\left( (L^{<0}\mathfrak{g})^*[1]\right)\otimes A\otimes \mathrm{Sym}^{\bullet}\!\!\left(\left(V(\CK)/V(\CO)\right)^*[2]\right).
\end{equation}
We comment that this should be understood as a projective system of resolutions. To make this a DG resolution, we need to include the usual Koszul differential $d_1$ coming from the pair $\widehat{L^{<0}\mathfrak{g}_0}$ and $(L^{<0}\mathfrak{g})^*$, as well as $d_2$ coming from the pair $ V(\CK)/V(\CO)$ and $(V(\CK)/V(\CO))^*$. However, these are not enough, since $\{D,d_2\}\ne 0$. One can in fact compute this commutator explicitly: let $u_{i,n}^*$ be the linear function on $V(\CK)/V(\CO)$ corresponding to $w_{i,n}^*$ in the Koszul resolution, then:
\begin{equation}
\{D,d_2\}u_{i,n}^*= Dw_{i,n}^*=\sum_{j, m+k=n}\rho_{i,m}^j\otimes v_{j,k}^*+\sum_{j_1,j_2, m_1+m_2+k=n} \rho_{i,m_1}^{j_1}\rho_{j_1, m_2}^{j_2}\otimes v_{j_2,k}^*+\cdots.
\end{equation}
And this commutator acts trivially on other generators of the pro-algebra. To make this into a DG resolution, we need to include another differential $\tilde{D}$: let $\epsilon_{i,n}^j$ be the linear function on $L^{<0}\mathfrak{g}$ corresponding to $\rho_{i,n}^j$ in the Koszul resolution. Then we define $\tilde{D}$ by:
\begin{equation}
\tilde{D}u_{i,n}^*=\sum_{j, m+k=n}\epsilon_{i,m}^j\otimes v_{j,k}^*+\frac{1}{2}\sum_{j_1,j_2, m_1+m_2+k=n}\epsilon_{i,m_1}^{j_1}\rho_{j_1, m_2}^{j_2}\otimes v_{j_2,k}^*+\cdots.
\end{equation}
This differential is of course $\tilde{G}_\CO$ invariant, since $\epsilon_{i,n}^j$ transforms in the same way as $\rho_{i,n}^j$, and $u_{i,n}^*$ transforms in the same way as $w_{i,n}^*$.  After introducing this new differential, $\{D,d_2\}=\{\tilde{D},d_1\}$ and the combination of the four differentials will be a differential, and the above indeed becomes a projective system of free resolutions. Now if we take endomorphism with $\mathbb{C}[V(\CO)]$, we obtain the space:
\begin{equation}
\mathrm{Sym}^{\bullet}\!\!\left(L^{<0}\mathfrak{g}[-1]\right)\otimes\mathbb{C}[V(\CO)]\otimes \mathrm{Sym}^{\bullet}\!\!\left(V(\CK)/V(\CO)[-2]\right),
\end{equation}
and the only nonzero differential is that induced from $\tilde{D}$. Examining the definition of $\tilde{D}$, we find that the higher order terms all drop off, and the linear term maps $L^{<0}\mathfrak{g}$ to $\mathbb{C}[V(\CO)]\otimes \mathrm{Sym}^{\bullet}\!\!\left(V(\CK)/V(\CO)[-2]\right)$, and is identified with the differential induced by the moment map. 

Combining the above steps, we obtain the following:

\newtheorem{Identitywithmatter}[LemHequiv]{Theorem}

\begin{Identitywithmatter}\label{Extidentitymatter}
There is a quasi-isomorphism of DG-$\mathbb{C}^*$ modules:
\begin{equation}\label{eqidentitywithmatter}
\begin{aligned}
\mathrm{End}^{\mathbb{B}\mathbb{C}^*}\!\!( i_*(\CO_{\overline{V(\CO)}}))\cong \left[ \mathbb{C}[V(\CO)]\otimes \mathrm{Sym}^{\bullet}\!\!\left(V(\CK)/V(\CO)[-2]\right) \otimes \mathrm{Sym}^{\bullet}\!\! \left(\left(\mathfrak{g}(\CK)/\mathfrak{g}(\CO)\oplus (z\mathfrak{g}(\CO))^*\right)[-1]\right) \right]^G.
\end{aligned}
\end{equation}
If we shift the loop grading of $V(\CK)$ by $q^{1/2}$, the cohomological degree of $V(\CK)/V(\CO)$ to $0$, and the cohomological degree of $\mathfrak{g}(\CK)/\mathfrak{g}(\CO)$ to $-1$, then the cohomology of this space coincides with $\pi_0\mathcal{V}_{G,V}$ in equation \eqref{eqpoissonmatter}.
\end{Identitywithmatter}

\begin{proof}

After the grading shift, $\mathrm{Sym}^{\bullet}\!\!\left(V(\CK)/V(\CO)\right)$ can be identified with $\mathbb{C}[V^*(\CO)]$. Since $\pi_0\mathcal{V}_{G,V}$ restricts to taking the Lie group invariants, comparing equation \eqref{eqidentitywithmatter} with \eqref{eqpoissonmatter}, we conclude that $\pi_0\mathcal{V}_{G,V}\cong \mathrm{End}^{\mathbb{B}\mathbb{C}^*}( i_*(\CO_{\overline{V(\CO)}}))$. 

\end{proof}

\noindent\textbf{Remark.} The character of the space in equation \eqref{eqidentitywithmatter} is given by:
\begin{equation}
\frac{1}{|W|}\oint_T \frac{\mathrm{d}s}{2\pi i s}\prod\limits_{\alpha\text{ roots}}  (1-s^\alpha) \left[(q)_{\infty}^{2\mathrm{rank}(G)}\frac{\prod\limits_{\alpha\text{ roots}}(qs^\alpha; q)_{\infty}^2}{\prod\limits_{\beta\text{ weights of } N\oplus N^*} (-q^{1/2}s^{\beta},q)_\infty}\right].
\end{equation}
This is the Schur index for the gauge theory with matter as stated in \cite{oh2020poisson}. 

\section{The Insertion of Fundamental t'Hooft Lines in Pure $PSL(2)$}\label{5}

Using the technique developed in the previous sections, one can consider the space of local operators at the junction of two half-BPS Wilson-'t Hooft line operators. As stated in the introduction, these operators correspond to vector bundles on the reduced $\tilde{G}_\CO$ orbits of $\CBFN$. Among these, the 't Hooft line operators are certain line bundles on the  $\tilde{G}_\CO$ orbits, and are labelled by the dominant coweight of $G$. These Wilson-'t Hooft line operators are the perverse coherent sheaves appearing in the work of \cite{cautis2019cluster, Cautistoappear}. For a full dictionary of correspondences between line operators and coherent sheaves, see \cite{kapustin2006holomorphic, cautis2019cluster, Cautistoappear}. 

Given two line operators $L_1$ and $L_2$, the space of local operators at their adjunction is given by:
\begin{equation}
\mathrm{Ops}_{G,V}(L_1,L_2)=\mathrm{Hom}_{\mathrm{Coh}(\tilde{G}_\CO\setminus \CBFN)}^{\mathbb{B}\mathbb{C}^*}(L_1,L_2).
\end{equation}
The space $\mathrm{Ops}_{G,V}(L_1,L_2)$ should give rise to a module of the Poisson vertex algebra $\mathrm{Ops}_{G,V}=\pi_0\mathcal{V}_{G,V}$. Indeed, since $\mathbbm{1}*L_i\cong L_i$ for $i=1,2$, the Poisson algebra $\mathrm{Ops}_{G,V}$ acts on $\mathrm{Ops}_{G,V}(L_1,L_2)$ through convolution. By the work of \cite{butson2021equivariant}, this action is also compatible with the factorization structure. The structure of these spaces as $\mathrm{Ops}_{G,V}$ modules has not been carefully described in literature; however, the Euler character of these spaces $\chi_q\mathrm{Ops}_{G,V}(L_1,L_2)$ are computed in \cite{cordova2016infrared}. 

We will look at the simplest non-trivial example: the space of local operators at the junction of fundamental 't Hooft lines and basic dyonic Wilson-'t Hooft lines in pure $PSL(2)$ ($PSU(2)$ in physics notation) theory. The fundamental t'Hooft line here is the structure sheaf of the miniscule orbit $\Gr_{1/2}\cong \mathbb{P}^1$, corresponding to the minimal dominant coweight $\frac{1}{2}$ of $PSL(2)$. We wiil not try to identify these as representations of $\mathrm{Ops}_{G,V}$, but only compute the vector spaces and their indices. We will match the indices with the indices of \cite{cordova2016infrared}.

We will keep using the notation $\overline{X}$ for the quotient stack $X/\tilde{G}_\CO$. Consider now $L_1=L_2=\CO_{\overline{\Gr_{1/2}}}$, we would like to compute:
\begin{equation}
\mathrm{Ops}_{G,V}(L_1,L_2)=\mathrm{End}^{\mathbb{B}\mathbb{C}^*}(\CO_{\overline{\Gr_{1/2}}})
\end{equation}
Proposition \ref{minisculeorbit} reduces the computation of Endomorphism algebra to computing the global sections of an associated vector bundle. In the computations below, we will drop all the quotients by $\tilde{G}_\CO$ in order to simlify the notations, although all the discussions below is carried in the equivariant settings. Fix a fixed point $z^{1/2}$, let $\tilde{P}$ be the stabilizer. The module $V$ is given by: 
\begin{equation}
\mathfrak{g}(\mathcal{K})/\left(z^{1/2}\mathfrak{g}(\mathcal{O})z^{-1/2}+\mathfrak{g}(\mathcal{O})\right)
\end{equation} 
The computation then requires that we understand the associated vector bundle as a bundle over $\mathbb{P}^1$. 

Since $z^{1/2}\mathfrak{g}(\mathcal{O})z^{-1/2}+\mathfrak{g}(\mathcal{O})=H(\mathcal{O})\oplus E(\mathcal{O})\oplus z^{-1}F(\mathcal{O})$, the representation:
\begin{equation}
 \mathfrak{g}(\mathcal{K})/\left(z^{1/2}\mathfrak{g}(\mathcal{O})z^{-1/2}+\mathfrak{g}(\mathcal{O})\right)
\end{equation}
falls into the following exact sequence: 
\begin{equation}
0\to z^{-1}\mathfrak{b}\to \mathfrak{g}(\mathcal{K})/\left(z^{1/2}\mathfrak{g}(\mathcal{O})z^{-1/2}+\mathfrak{g}(\mathcal{O})\right)\to  \mathfrak{g}(\mathcal{K})/z^{-1}\mathfrak{g}(\CO)\to 0,
\end{equation}
where $\mathfrak{b}$ is the Lie algebra of $B\subseteq G$. This short exact sequence split as a representation of $B$, and since $G(\CO)/\tilde{P}=G/B$, we have an isomorphism of vector bundles:
\begin{equation}
\tilde{G}_\CO \times_{\tilde{P}} V \cong G\times_B z^{-1}\mathfrak{b}\bigoplus \CO_{\mathrm{Gr}_{1/2}}\otimes \mathfrak{g}(\mathcal{K})/z^{-1}\mathfrak{g}(\CO).
\end{equation}
Taking exterior power on both sides, we have an isomorphism:
\begin{equation}
\tilde{G}_\CO \times_{\tilde{P}} \mathrm{Sym}^{\bullet}\!\left( V[1]\right)\cong G\times_B  \mathrm{Sym}^{\bullet}\!(z^{-1}\mathfrak{b}[-1])\otimes    \mathrm{Sym}^{\bullet}\!\!\left(\mathfrak{g}(\mathcal{K})/z^{-1}\mathfrak{g}(\CO)[-1]\right).
\end{equation}
The global section of $G\times_B  \mathrm{Sym}^{\bullet}\!(z^{-1}\mathfrak{b}[-1])$ can be computed easily:
\begin{equation}
\begin{aligned}
& H^*(G\times_B  \mathrm{Sym}^{0}\!(z^{-1}\mathfrak{b}[-1]))=\mathbb{C}[0],\\
&q^{-1} H^*(G\times_B  \mathrm{Sym}^{1}\!(z^{-1}\mathfrak{b}[-1]))=\mathbb{C}[-1]\oplus  \mathfrak{g}[-1],\\
&q^{-2}H^*(G\times_B  \mathrm{Sym}^{2}\!(z^{-1}\mathfrak{b}[-1]))=\mathfrak{g}[-2].
\end{aligned}
\end{equation}
The index of $H^*(G\times_B  \mathrm{Sym}^{\bullet}\!(z^{-1}\mathfrak{b}[-1]))$ is equal to $(1-q)(1-q-qs^2-qs^{-2})$.  The index of $\mathrm{Sym}^{\bullet}\!\!\left(\mathfrak{g}(\mathcal{K})/z^{-1}\mathfrak{g}(\CO)[-1]\right)$ is given by:
\begin{equation}
\frac{(q)_\infty (q:qs^2)_\infty(q:qs^{-2})_\infty}{(1-q)(1-qs^2)(1-qs^{-2})},
\end{equation}
so the index of $\mathrm{End}^{\mathrm{Rep}(\tilde{G}_\CO)}(\CO_{\overline{\Gr_{1/2}}})$ is given by:
\begin{equation}
 \frac{1-q-qs^2-qs^{-2}}{(1-qs^2)(1-qs^{-2})} (q)_\infty (q:qs^2)_\infty(q:qs^{-2})_\infty.
\end{equation}
In conclusion:

\newtheorem{Thmfundamental}{Theorem}[section]

\begin{Thmfundamental}\label{Thmfundamental}
The space $\mathrm{End}^{\mathbb{B}\mathbb{C}^*}(\CO_{\overline{\Gr_{1/2}}})$ is quasi-isomorphic to:
\begin{equation}\label{localopthooft}
\left[\left(\mathbb{C}\oplus q\mathbb{C}[-1]\oplus q\mathfrak{g}[-1]\oplus q^2\mathfrak{g}[-2]\right)\otimes \mathrm{Sym}^{\bullet}\!\!\left(\left(\mathfrak{g}(\mathcal{K})/z^{-1}\mathfrak{g}(\CO)\oplus (z\mathfrak{g}(\CO))^*\right)[-1]\right)\right]^G.
\end{equation}
The character of this space is:
\begin{equation}
\frac{1}{2}(q)_\infty^2\oint_s\mathrm{d}s (1-s^2)(1-s^{-2})\frac{(1-q-qs^2-qs^{-2})}{(1-qs^2)(1-qs^{-2})}(q:qs^2)_\infty^2(q:qs^{-2})_\infty^2.
\end{equation}
\end{Thmfundamental}

We expect \eqref{localopthooft} to be the space of local operators supported at a point on a single straight 't Hooft line. The resulting character does not match the index obtained in \cite{cordova2016infrared}. The computation that matches their result is the following: consider now $L_1=\CO_{\overline{\Gr_{1/2}}}$ the fundamental 't Hooft line, and $L_2=\Omega_{\overline{\Gr_{1/2}}}$ a dyonic Wilson-'t Hooft line, where $\Omega_{\overline{\Gr_{1/2}}}$ is the dualizing sheaf of $\Gr_{1/2}$. Physically, $\Omega_{\overline{\Gr_{1/2}}}$ corresponds to the dyonic Wilson-'t Hooft line with fundamental magnetic charge and $+1$ electric~charge. Consider the junction:
\begin{equation}
\mathrm{Ops}_{G,V}(L_1,L_2)=\mathrm{Hom}^{\mathbb{B}\mathbb{C}^*}(\CO_{\overline{\Gr_{1/2}}},\Omega_{\overline{\Gr_{1/2}}}),
\end{equation}
In this case, Proposition \ref{minisculeorbit} still applies, with the associated bundle of $V$ twisted by the canonical sheaf of $\Gr_{1/2}=G/B$. The cohomology of the twisted $G\times_B z^{-1}\mathfrak{b}\otimes O(-2)$ is:
\begin{equation}
\begin{aligned}
& H^*(G\times_B \mathrm{Sym}^0(z^{-1}\mathfrak{b}\otimes O(-2)[-1]))=\mathbb{C}[-1],\\
&q^{-1} H^*(G\times_B \mathrm{Sym}^1(z^{-1}\mathfrak{b}\otimes O(-2)[-1]))=\mathbb{C}[-1]\oplus  \mathbb{C}[-2],\\
&q^{-2}H^*(G\times_B \mathrm{Sym}^2(z^{-1}\mathfrak{b}\otimes O(-2)[-1]))=\mathbb{C}[-2].
\end{aligned}
\end{equation}
This space has index $-(1-q)(1+q)$. The contribution from $\mathfrak{g}(\CK)/z^{-1}\mathfrak{g}(\CO)$ remains the same. Hence the space $\mathrm{Hom}^{\mathbb{B}\mathbb{C}^*}(\CO_{\overline{\Gr_{1/2}}},\Omega_{\overline{\Gr_{1/2}}})$ has index:
\begin{equation}
-\frac{1}{2}(q)_\infty^2\oint_s\mathrm{d}s (1-s^2)(1-s^{-2})\frac{(1+q)}{(1-qs^2)(1-qs^{-2})}(q:qs^2)_\infty^2(q:qs^{-2})_\infty^2,
\end{equation}
Shifting by $q^{1/2}$, this exactly matches the formula in \cite{cordova2016infrared}. In this paper, the authors are implicitly using a Serre functor to rotate the line operators, or in other words, they took the dual of the line operators in the monoidal category. This subtle operation was described explicitly in \cite{cautis2019cluster}, and an important feature is that the left dual of a line operator is \textbf{not} necessarily equivalent to the right dual, unless the theory is superconformal, cf. the end of Section \ref{3}. As explained in \cite{cautis2019cluster}, the difference between left dual and right dual is due to the fact that the dualizing sheaf $\Omega$ of $\tilde{G}_\CO\setminus \Gr_G$ is not isomorphic to its involution $s^*\Omega$, where the involution $s: \tilde{G}_\CO\setminus \Gr_G\to \tilde{G}_\CO\setminus \Gr_G$ is defined by $s([g])=[g^{-1}]$. Thus, after rotating by $2\pi$, a line operator receives contribution from the dualizing sheaf of $\Gr_G$. In the case of the structure sheaf of a miniscule orbit, this is simply the dualizing sheaf of the orbit. This explains the presence of $\Omega_{\overline{\Gr_{1/2}}}$ in the above formula. 

\newpage

\bibliographystyle{plain}

\bibliography{formal}

\end{document}